\documentclass{lmcs}
\pdfoutput=1

\usepackage{lastpage}
\lmcsdoi{20}{4}{2}
\lmcsheading{}{\pageref{LastPage}}{}{}%
{Jul.~19,~2022}{Oct.~07,~2024}{}

\usepackage{amsfonts,qtree,stmaryrd,textcomp}
\usepackage[bbgreekl]{mathbbol}
\usepackage{pifont,amssymb,amsmath,amsthm,tabularx,graphicx}
\usepackage{tree-dvips,pictexwd,dcpic}
\usepackage{bbm}
\usepackage{bm}
\usepackage[T1]{fontenc}
\usepackage[latin1]{inputenc}
\usepackage{fancyhdr}
\usepackage{graphicx}
\usepackage{enumitem}

\usepackage[all]{xy}
\usepackage{epigraph}
\RequirePackage{bussproofs}

\usepackage{mathabx}
\usepackage{scalerel,stackengine}

\usepackage{framed}
\usepackage{mathtools}
\usepackage{url}
\usepackage{hyperref}
\usepackage{proof}
\usepackage{xspace}
\usepackage{enumitem}
\usepackage{listings}
\usepackage{wrapfig}
\usepackage{tikz}
\usepackage{tikz-cd}
\usepackage{tikz-3dplot}

\usepackage{relsize}

\usetikzlibrary{positioning}
\usetikzlibrary{shapes}
\usetikzlibrary{arrows}
\usetikzlibrary{calc}
\usetikzlibrary{decorations.pathmorphing}

\usetikzlibrary{decorations.pathreplacing}
\tikzset{axis/.style={&lt;-&gt;}}

\stackMath
\newcommand\reallywidehat[1]{%
\savestack{\tmpbox}{\stretchto{%
  \scaleto{%
    \scalerel*[\widthof{\ensuremath{#1}}]{\kern-.6pt\bigwedge\kern-.6pt}%
    {\rule[-\textheight/2]{1ex}{\textheight}}
  }{\textheight}%
}{0.5ex}}%
\stackon[1pt]{#1}{\tmpbox}%
}

\definecolor{MyBlue}{rgb}{0.05, 0.25, 0.65}
\definecolor{MyRed}{rgb}{0.90, 0.05, 0.05}
\definecolor{MyGreen}{rgb}{0.05, 0.90, 0.05}

\newcommand{\B}{\boldsymbol}
\newcommand{\C}[1]{\mathcal{#1}}
\newcommand{\D}[1]{\mathbb{#1}}
\newcommand{\F}[1]{\mathfrak{#1}}

\newcommand{\Nat}{{\mathbb N}}
\newcommand{\Real}{{\mathbb R}}

\newcommand{\id}{\mathrm{id}}

\newcommand{\BISH}{\mathrm{BISH}}

\newcommand{\CST}{\mathrm{CST}}

\newcommand{\dom}{\mathrm{dom}}

\newcommand{\TOT}{\Leftrightarrow}

\newcommand{\To}{\Rightarrow}

\newcommand{\sto}{\rightsquigarrow}

\newcommand{\MLTT}{\mathrm{MLTT}}
\newcommand{\CZF}{\mathrm{CZF}}

\newcommand{\Fam}{\textnormal{\texttt{Fam}}}

\newcommand{\pr}{\textnormal{\texttt{pr}}}

\newcommand{\BST}{\mathrm{BST}}

\newcommand{\Disj}{\B \rrbracket \B \llbracket}

\newcommand{\Set}{\mathrm{\mathbf{Set}}}

\newcommand{\se}{\mathrm{se}}

\newcommand{\eto}{\hookrightarrow}

\newcommand{\BMT}{\mathrm{BMT}}

\newcommand{\BCMT}{\mathrm{BCMT}}

\newcommand{\PBCMT}{\textnormal{\texttt{PBCMT}}}

\newcommand{\BCMS}{\textnormal{\texttt{BCMS}}}

\newcommand{\PIS}{\mathrm{PIS}}

\newcommand{\PMS}{\mathrm{PMS}}

\newcommand{\BCIS}{\mathrm{BCIS}}

\newcommand{\supp}{\mathrm{supp}}

\newcommand{\Supp}{\mathrm{Supp}}

\newcommand{\Simple}{\textnormal{\texttt{Simple}}}

\newcommand{\disjrep}{\textnormal{\texttt{disjrep}}}

\keywords{constructive measure theory, measure spaces, integration spaces, integrable functions}

\begin{document}

\title[Pre-measure and pre-integration spaces in BCMT]{Pre-measure spaces and pre-integration spaces in predicative Bishop-Cheng measure theory}

\author[I.~Petrakis]{Iosif Petrakis\lmcsorcid{0000-0002-4121-7455}}[a]
\author[M.~Zeuner]{Max Zeuner\lmcsorcid{0000-0003-3092-8144}}[b]

\address{University of Verona, Department of Computer Science}
\email{iosif.petrakis@univr.it}
\thanks{The research of the first author
was partially supported by LMUexcellent, funded by the Federal
Ministry of Education and Research (BMBF) and the Free State of Bavaria under the
Excellence Strategy of the Federal Government and the L\"ander.}

\address{Stockholm University, Mathematics Department}
\email{zeuner@math.su.se}

\begin{abstract}
\noindent
Bishop's measure theory (BMT), introduced in~\cite{Bi67}, is an
abstraction of the measure theory of a locally compact metric space
$X$, and the use of an informal notion of a set-indexed family of
complemented subsets is crucial to its predicative character. The more
general Bishop-Cheng measure theory (BCMT), introduced in~\cite{BC72}
and expanded in~\cite{BB85}, is a constructive version of the
classical Daniell approach to measure and integration, and highly
impredicative, as many of its fundamental notions, such as the
integration space of $p$-integrable functions $L^p$, rely on
quantification over proper classes (from the constructive point of
view).  In this paper we introduce the notions of a pre-measure and
pre-integration space, a predicative variation of the Bishop-Cheng
notion of a measure space and of an integration space, respectively.
Working within Bishop Set Theory $(\BST)$, elaborated in~\cite{Pe20},
and using the theory of set-indexed families of complemented subsets
and set-indexed families of real-valued partial functions within
$\BST$, we apply the implicit, predicative spirit of BMT to BCMT. As a
first example, we present the pre-measure space of complemented
detachable subsets of a set $X$ with the Dirac-measure, concentrated
at a single point. Furthermore, we translate in our predicative
framework the non-trivial, Bishop-Cheng construction of an integration
space from a given measure space, showing that a pre-measure space
induces the pre-integration space of simple functions associated to
it. Finally, a predicative construction of the canonically integrable
functions $L^1$, as the completion of an integration space, is
included.
\end{abstract}

\maketitle

\section{Introduction}
\label{sec: intro}

In the most popular approach to classical measure theory, see
e.g.,~\cite{Ha74}, integration is defined through measure.  Starting
from a measure space $(X, \C A, \mu)$, one defines simple and
measurable functions, the latter through the Borel sets in $\Real$. As
a positive measurable function is the limit of an increasing sequence
of positive, simple functions, the obviously defined integral of a
simple function is extended to the integral of a positive, measurable
function.  The integral of a measurable function $f \colon X \to
\Real$ is then defined through the integrals of the positive,
measurable functions $f_+$ and $f_{-}$. The highly non-constructive
standard approach can be roughly characterised as an approach ``from
sets to functions''.

In the \textit{Daniell approach} to classical measure theory, see
e.g.,~\cite{Lo53, Ta73}, measure is defined through integration. It
was introduced by Daniell~\cite{Da18}, it was taken further by
Weil~\cite{We40}, Kolmogoroff~\cite{Ko48}, Stone~\cite{St48},
Carath\'eodory~\cite{Ca56}, and Segal~\cite{Se54, Se65}, and it is
incorporated in Bourbaki~\cite{Bo04}.  The starting point of the
Daniell approach is the notion of \textit{Daniell space} $\big(X, L,
\int\big)$, where $L$ is a Riesz space of real-valued functions on $X$
and $\int \colon L \to \Real$ is a positive, linear functional that
satisfies the Daniell property, a certain continuity condition.  Using
the (non-constructive) Bolzano-Weierstrass theorem, one extends $L$ to
$L^+$, which is the set of functions $f \colon X \to \overline{\Real}$
that are limits of increasing sequences in $L$, and $\int$ is extended
to $\int^+ \colon L^+ \to \overline{\Real}$ accordingly. The upper
$\overline{\int} f$ and lower integral $\underline{\int} f $ of a
function $f \colon X \to \overline{\Real}$ are defined through the
(non-constructive) completeness axiom of real numbers, and $f$ is
integrable, or an element of $L^1$, if
$\underline{\int} f = \overline{\int}f\in\Real$.
A function $f \colon X \to [0, + \infty]$
is called \textit{measurable}, if it can be approximated appropriately
by integrable functions, and a subset $A$ of $X$ is measurable, if its
characteristic function $\chi_A$ is measurable, while $A$ is
\textit{integrable}, if $\chi_A \in L^1$. If $A$ is integrable, a
measure function $A \mapsto \mu(A)$ is defined through the integral of
$\chi_A$. A clear advantage of this approach is that ``certain
properties of the integral already follow from the integrals of the
nice functions, which are easier to handle than arbitrary integrable
functions''~\cite{Wi90}. The Daniell approach can be roughly
characterised as an approach ``from functions to sets''.

As functions are more appropriate to constructive study than sets,
Bishop followed the Daniell approach both in~\cite{Bi67}, and, in a
different and more uniform way, in~\cite{BC72, BB85}.  We call
\textit{Bishop measure theory} $(\BMT)$ the measure theory developed
by Bishop in~\cite{Bi67}.  Although the integration theory of locally
compact metric spaces within $\BMT$ follows the Daniell approach, the
treatment of abstract measures within $\BMT$ follows the more popular
approach to classical measure theory. As the Borel sets are defined
inductively in $\BMT$, the set theory required for it must accommodate
inductive definitions with rules of countably many premises.

The more general theory of measure introduced in~\cite{BC72}, and
significantly extended in~\cite{BB85}, is what we call
\textit{Bishop-Cheng measure theory} $(\BCMT)$, which makes no use of
(inductively defined) Borel sets, and hence it is based on a set
theory without inductive definitions.  Following the tradition of the
Daniell approach, Bishop and Cheng consider first the integral on a
certain set $L$ of given functions, then extend it to the larger set
of functions $L^1$, and define the measure at a later stage. Although
complemented subsets\footnote{These are pair of subsets that are
disjoint in a positive and strong way. Their use in $\BMT$ and $\BCMT$
is crucial to avoid many negatively defined concepts from their
classical counterparts.} are first-class citizens both in $\BMT$ and
in $\BCMT$, their set-indexed families are not employed in $\BCMT$.
What we call here the \textit{Bishop-Cheng integration space} is the
constructive analogue to Daniell space\footnote{In~\cite{Pe22c} it is
explained why the notion of a Bishop-Cheng integration space is a
natural, constructive counterpart to the classical notion of Daniell
space.}  that captures all basic examples of the classical Daniell
theory.  The broadness of results within $\BCMT$ presented
in~\cite{BB85} and in several related publications is
striking. Numerous applications of Bishop-Cheng measure theory to
probability theory and to the theory of stochastic processes are found
in the older work of Chan~\cite{Ch72}-\cite{Ch75}, and especially in
his recent monograph~\cite{Ch21}. The generality of $\BCMT$ though, is
due to the use of impredicative definitions, which hinder the
extraction of efficient computations from proofs.

If $\F F^{se}(X)$ is the totality of strongly extensional,
real-valued, partial functions $f$ on a set with a given inequality
$(X, =_X, \neq_X)$, the set of integrable functions $L^1$ is defined
in $\BCMT$ by the separation scheme as follows:
$$L^1 := \big\{f \in \F F^{se}(X) \mid f \ \mbox{is integrable}\big\}.$$
The membership-condition of the totality $\F F^{se}(X)$ involves
quantification over the universe of sets, since a partial, real valued
function is by definition a set $A$ together with an embedding (or
injection) $i_A$ of $A$ into $X$ and a function $f \colon A \to
\Real$. Hence, $\F F^{se}(X)$ is a proper class, and the separation
scheme on a proper class does not define a set. Thus, from a
predicative point of view, the Bishop-Cheng definition of $L^1$ does
not determine a set.  As this impredicativity of $L^1$ is ``dense'' in
$\BCMT$, the original approach of Bishop and Cheng, as a whole, cannot
express successfully the computational content of measure
theory. Exactly this computational deficiency of $\BCMT$ is also
recognised by Spitters in~\cite{Sp06}.

Already in the definition of a Bishop-Cheng integration space a
similar problem arises. Namely, the integral is supposed to be defined
on a subset $L$ of the proper class $\F F^{se}(X)$, without specifying
though, how such a subset can be defined i.e., how a subclass of $\F
F^{se}(X)$ can be considered to be a set. It seems that both
in~\cite{BC72} and in~\cite{BB85} the totality $\F F^{se}(X)$ is taken
to be a set. This fundamental impredicativity built in $\BCMT$
directed the subsequent constructive studies of measure theory to
different directions\footnote{Outside Bishop's constructivism there
are various approaches to measure theory.  The theory of
measure~\cite{He56} within Brouwer's intuitionism contradicts the
classical theory, while measure theory~\cite{Ed09} within the
computability framework of Type-2 Theory of Effectivity is based on
classical logic. Measure theory~\cite{Sa68},~\cite{BD91} within
Russian constructivism employs Markov's principle of unbounded
search. In intuitionisitic Martin-L\"of type theory
$(\MLTT)$~\cite{ML98} the interest lies mainly in probabilistic
programming~\cite{BAV12}, while in homotopy type theory~\cite{HoTT13}
univalent techniques, such as higher inductive types, are applied to
probabilistic programming too~\cite{FS16}.}.

Coquand, Palmgren, and later Spitters, also acknowledged that $\BCMT$
does not facilitate the extraction of efficient
computations. According to Spitters~\cite{Sp06}, it is unlikely that
$\BCMT$ ``will be useful when viewing Bishop-style mathematics as a
high-level programming language''. As a result, the search of the
computational content of measure theory in constructive mathematics
was shifted from the Bishop-Cheng theory to more abstract, algebraic,
or point-free approaches (see the work of Coquand, Palmgren and
Spitters, in~\cite{CP02},~\cite{Sp06} and~\cite{CS09}).  However, in
terms of applications\footnote{The applications to probability theory
were difficult to explore in Spitters' approach and postponed in the
approach of Coquand and Palmgren.  Recently, a decisive step towards a
point-free treatment of measure theory has been taken by
Simpson~\cite{Si12}.  Simpson advocates however in a classical
framework, that sublocales, rather than subspaces, be vital. This
conceptual move allows one to even circumvent some of the constraints
at the outset of measure theory.  A fairly constructive development,
which had been kept by Simpson for future work, has been proposed by
Ciraulo~\cite{Ci22}.  However, Ciraulo invokes the principle of
countable choice. Interestingly, he is also concerned with the status
of complemented sets from a point-free perspective.}, these approaches
attain neither the range nor the broadness of $\BCMT$.

Already in $\BMT$ though, Bishop avoided impredicativities by using
(two) set-indexed families of complemented subsets in his definition
of a measure space, in order to quantify over the index-sets
only. Discussing in~\cite{Bi70}, p.~67, the exact definition of a
measure space in $\BMT$ within his formal system $\Sigma$, he writes
the following:
\begin{quote}
To formalize in $\Sigma$ the notion of an abstract measure space,
definition 1 of chapter 7 of~\cite{Bi67} must be rewritten as
follows. A \textit{measure space} is a family $\C M \equiv \{A_t\}_{t
  \in T}$ of complemented subsets of a set $X$ $\ldots$, a map $\mu :
T \to \Real^{0+}$, and an additional structure $\ldots$ . If $s$ and
$t$ are in $T$, there exists an element $s \mathsmaller{\vee} t$ of
$T$ such that $A_{s \mathsmaller{\vee} t} < A_s \cup A_t$.  Similarly,
there exist operations $\mathsmaller{\wedge}$ and $\mathsmaller{\sim}$
on $T$, corresponding to the set theoretic operations $\cap$ and $-$.
The usual algebraic axioms are assumed, such as
$\mathsmaller{\sim}(s \mathsmaller{\vee} \ t) = \mathsmaller{\sim} s\ \mathsmaller{\wedge} \ \mathsmaller{\sim} t$. $\ldots$
Considerations such as the above indicate that essentially all of the
material in~\cite{Bi67}, appropriately modified, can be comfortably
formalised in $\Sigma$.
\end{quote}
\noindent
This indexisation method, roughly sketched in~\cite{Bi67}, is
elaborated within Bishop Set Theory $(\BST)$ in~\cite{Pe20}.  Based on
this, we present here the first crucial steps to a predicative
reconstruction $(\PBCMT)$ of $\BCMT$. Following Bishop's explanations
in~\cite{Bi70}, we replace a totality of strongly extensional,
real-valued, partial functions $L$ in the original definition of a
Bishop-Cheng integration space by a set-indexed family $\Lambda$ of
such partial functions.  Applying tools and results from~\cite{Pe20},
we recover the concept of an integration space in an indexised
form. The predicative advantage of the indexisation method within
$\PBCMT$ is that crucial quantifications are over an index-set and not
proper classes. Following~\cite{Pe20}, we elaborate the concept of a
pre-integration space in which the index-set $I$ is equipped with all
necessary operations so that a pre-integral $\int$ can be defined on
$I$. A pre-integration space induces a predicative integration space,
the integral $\int^*$ of which on the partial function $f_i$ is given,
for every $i \in I$, by
$$\int^* f_i := \int i.$$

We provide a predicative treatment of $L^1$ by considering only the
canonically integrable functions\footnote{This terminology is
introduced by Spitters in~\cite{Sp02}.}  of a given pre-integration
space. Our main result is that the set-indexed family of canonically
integrable functions admits the structure of a pre-integration space
(Theorem~\ref{thm: pre-2.18}), which is an appropriate completion of
the original pre-integration space (Theorem~\ref{thm: pre-2.16}).  The
theory developed in~\cite{Pe20} together with careful arguments that
avoid the use of the class of full sets and countable choice
(see~\cite{Ri01, Sc04} for a critique to the use of countable choice
in Bishop-style constructive mathematics, also known as $\BISH$)
helped us prove a constructive and predicative version of Lebesgue's
series theorem (Theorem~\ref{thm: Leb}), which is crucial to the proof
of our main result.  A predicative definition of $L^1$ ensures that
all concepts defined through quantification over $L^1$ in $\BCMT$ are
also predicative. For example, quantification over $L^1$ is used in
the Bishop-Cheng definition of a full set, which is a constructive
counterpart to the complement of a null set in classical measure
theory. This predicative treatment of $L^1$ is the first, clear
indication that the computational content of measure theory can be
grasped by the predicative reconstruction $\PBCMT$ of the original
$\BCMT$.

\section{Overview of this paper}
\label{sec: overview}

We structure this paper as follows:

\begin{itemize}
\item In section~\ref{sec: complsub} we describe the connection
  between complemented subsets and boolean-valued partial functions,
  which explains the crucial role of partial functions in $\BCMT$. The
  constructive way to employ the passage from functions to sets in the
  classical Daniell approach through the use of characteristic
  functions of subsets, is to work with complemented subsets and their
  (partial) characteristic functions.

\item In section~\ref{sec: famsubsets} we describe the basic
  properties of set-indexed families of subsets of a given set $X$.
  We discuss the set-character within $\BST$ of the totality of
  families of subsets of $X$ indexed by some set $I$, which will be
  relevant to our presentation of a pre-measure space.

\item In section~\ref{sec: fampartial} we define within $\BST$ the
  notions of a family of partial functions and of a family of
  complemented subsets indexed by some set $I$. These
  function-theoretic concepts will be used in $\PBCMT$ instead of the
  abstract sets of partial functions and of subsets, respectively,
  that are considered in $\BCMT$.

\item In section~\ref{sec: premeasures} we introduce the notion of a
  pre-measure space as a predicative counterpart to the notion of
  Bishop-Cheng measure space in $\BCMT$. The pre-measure space of
  complemented detachable subsets of a set $X$ with the Dirac-measure
  concentrated at a single point is studied.

\item In section~\ref{sec: realpartial} we include the facts on
  real-valued, partial functions that are necessary to the definition
  of a pre-integration space within $\BST$ (Definition~\ref{def: preintspace}).

\item In section~\ref{sec: preint} we introduce the notion of a
  pre-integration space as a predicative counterpart to the notion of
  an integration space in $\BCMT$. We also briefly describe the
  pre-integration space $\big(X, I, \int {} d\mu\big)$, where $X$ is a
  locally compact metric space $X$ with a so-called modulus of local
  compactness, $I$ is the set of functions with compact support on
  $X$, and the integral $\int {f} d\mu$ of $f \in I$ is the measure
  $\mu(f)$, where $\mu$ is a positive measure on $X$
  (Theorem~\ref{thm: intlcms}).

\item In section~\ref{sec: simple} we construct the pre-integration
  space of simple functions from a given pre-measure space
  (Theorem~\ref{thm: pre-10.10}). This is a predicative translation
  within $\BST$ of the construction of a Bishop-Cheng integration
  space from the simple functions of a measure space (Theorem 10.10
  in~\cite{BB85}). Although we follow the corresponding construction
  in section 10 of chapter 6 in~\cite{BB85} closely, our approach
  allows us to not only work completely predicatively, but also to
  carry out all proofs avoiding the axiom of countable choice.

\item In section~\ref{sec: L1} we first present the canonically
  integrable functions explicitly as a family of partial functions, in
  order to avoid the impredicativities of the original Bishop-Cheng
  definition of $L^1$. Based on a predicative version of Lebesgue's
  series theorem (Theorem~\ref{thm: Leb}), we then show that this
  family admits the structure of a pre-integration space
  (Theorem~\ref{thm: pre-2.18}) and explain in what sense it can be
  seen as the completion of our original pre-integration space
  (Theorem~\ref{thm: pre-2.16}).

\item In section~\ref{sec: concl} we list some question for future
  work stemming from the material presented here.

\end{itemize}
\noindent
We work within $\BST$, which behaves as a high-level programming
language.  For all notions and results of Bishop set theory that are
used here without definition or proof we refer to~\cite{Pe21}, in this
journal\footnote{In~\cite{Pe21} the theory of spectra of Bishop spaces
  (see~\cite{Pe15}-\cite{Pe19b} and~\cite{Pe20b}-\cite{Pe22a}) is
  developed within $\BST$.},
and to~\cite{Pe20, Pe22b}.  For all notions
and results of constructive real analysis that are used here without
definition or proof we refer to~\cite{BB85}. The type-theoretic
interpretation of Bishop's set theory into the theory of setoids (see
especially the work of Palmgren~\cite{Pa05}-\cite{PW14}) has become
nowadays the standard way to understand Bishop sets\footnote{For an
  analysis of the relation between intensional $\MLTT$ and Bishop's
  theory of sets see~\cite{Pe20}, Chapter 1.}.
Other suitable, yet
different, formal systems for $\BISH$ are Myhill's Constructive Set
Theory $(\CST)$, introduced in~\cite{My75}, and Aczel's system $\CZF$
(see~\cite{AR10}).

\section{Partial functions and complemented subsets}
\label{sec: complsub}

Bishop set theory $(\BST)$, elaborated in~\cite{Pe20}, is an informal,
constructive theory of totalities and assignment routines that serves
as a ``completion'' of Bishop's original theory of sets in~\cite{Bi67,
  BB85}.  Its first aim is to fill in the ``gaps'', or highlight the
fundamental notions that were suppressed by Bishop in his account of
the set theory underlying Bishop-style constructive mathematics
$\BISH$.  Its second aim is to serve as an intermediate step between
Bishop's theory of sets and an \textit{adequate} and \textit{faithful}
formalisation of $\BISH$ in Feferman's sense~\cite{Fe79}. To assure
faithfulness, we use concepts or principles that appear, explicitly or
implicitly, in $\BISH$.  $\BST$ ``completes'' Bishop's theory of sets
in the following ways. It uses explicitly a universe of (predicative)
sets $\D V_0$, which is a proper class. It separates clearly sets from
proper classes. Dependent operations, which were barely mentioned
in~\cite{Bi67, BB85}, are first-class citizens in $\BST$. An
elaborated theory of set-indexed families of sets is included in
$\BST$. As an introduction to the basic concepts of $\BST$ is included
in~\cite{Pe21}, in this journal, and in~\cite{Pe20, Pe22b}, we refer
the reader to these sources for all basic concepts and results within
$\BST$ that are mentioned here without further explanation or
proof. Next we present some basic properties of partial functions and
complemented subsets within $\BST$, which are necessary to the rest of
this paper. A \textit{subset} of a set $X$ is a pair $(A, i_A)$, where
$(A, =_A)$ is a set and $i_A \colon A \eto X$ is an embedding i.e.,
$i_A(a) =_X i_A(a{'}) \To a =_A a{'}$, for every $a, a{'} \in A$. The
intersection of two subsets is given by the corresponding pullback,
and their union is defined in~\cite{Bi67}, p.~64. We denote the set of
functions from $A$ to $X$ by $\D F(A, X)$.

\begin{defi}\label{def: partialfunction}
Let $X, Y$ be sets.  A partial function\index{partial function} from
$X$ to $Y$ is a triplet $\B f_A := (A, i_A, f_A)$\index{$(A, i_A^X, f_A^Y)$},
where $(A, i_A) \subseteq X$, and $f_A \in \D F(A, Y)$.
We call $f_A$ total, if $\dom(\B f_A) := A =_{\C P(X)} X$.  Let $\B
f_A \leq \B f_B$\index{$f_A^Y \leq f_B^Y$}, if there is an embedding
$e_{AB} \colon A \eto B$ such that the following triangles commute
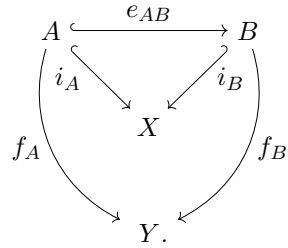
\begin{center}
\resizebox{4cm}{!}{%
\begin{tikzpicture}

\node (E) at (0,0) {$A$};
\node[right=of E] (B) {};
\node[right=of B] (F) {$B$};
\node[below=of B] (A) {$X$};
\node[below=of A] (C) {$ \  Y$.};

\draw[right hook->] (E)--(F) node [midway,above] {$e_{AB}$};
\draw[->,bend right=40] (E) to node [midway,left] {$f_A$} (C);
\draw[->,bend left=40] (F) to node [midway,right] {$f_B$} (C);
\draw[right hook->] (E)--(A) node [midway,left] {$i_A \ $};
\draw[left hook->] (F)--(A) node [midway,right] {$\ i_B  $};

\end{tikzpicture}
}
\end{center}
In this case we write $e_{AB} \colon \B f_A \leq \B f_B$.
The partial function space $\F F(X,Y)$ is
equipped with the equality
$\B f_A =_{\F F(X,Y)} \B f_B :\TOT \B f_A \leq \B f_B \ \& \ \B f_B \leq \B f_A$.
If $X, Y$ are equipped with inequalities $\neq_X, \neq_Y$, respectively, let
$\F F^{\se}(X, Y)$ be the totality\footnote{
  As the membership condition for $\F F(X,Y)$
  requires quantification over the universe of sets $\D V_0$,
  the totalities $\F F(X,Y)$ and $\F F^{\se}(X, Y)$ are proper classes.}
of strongly extensional elements of $\F F(X, Y)$.
\end{defi}

\begin{defi}\label{def: oppartial}
If $\D 2 := \{0, 1\}$ and $\B f = (A, i_A, f_A), \B g = (B, i_B, g_B) \in \F F(X, \D 2)$, let
$\B f \vee \B g := \max\{\B f, \B g\} = (A \cap B, i_{A \cap B}, f_A \vee g_B)$, $\B f \cdot \B g := \B f \wedge \B g := \min\{\B f, \B g\} := (A \cap B, i_{A \cap B}, f_A \wedge g_B)$,
$\sim \B f := 1 - \B f : = (A, i_A, 1 - f_A)$, and $\B f \sim \B g := \B f \wedge (\sim\B g)$,
where $1$ also denotes the constant function on $A$ with value $1$.
\end{defi}

An inequality on a set $X$ induces a positively defined notion of disjointness of subsets of $X$,
which in turn induces the notion of a complemented subset of $X$. In this way the negatively defined
notion of the set-theoretic complement of a subset is avoided.

\begin{defi}\label{def: apartsubsets}
Let $(X, =_X, \neq_X)$ be a set with inequality, and $(A, i_A), (B,
i_B) \subseteq X$. We say that $A$ and $B$ are disjoint with respect
to $\neq_X$\index{disjoint subsets}, in symbols $A \Disj B$\index{$A
  \Disj_{\mathsmaller{\neq}} B$}, if $\forall_{a \in A}\forall_{b \in
  B}\big(i_A(a) \neq_X i_B(b) \big)$.  A \textit{complemented
  subset}\index{complemented subset} of $X$ is a pair $\B A := (A^1,
A^0)$\index{$\B A := (A^1, A^0)$}, where $(A^1, i_{A^1}), (A^0,
i_{A^0}) \subseteq X$, such that $A^1 \Disj A^0$.  The
\textit{characteristic function} of $\B A$ is the operation\footnote{A non-dependent assignment routine
  $f:A\sto B$, where $A$ and $B$ are sets, is called an \textit{operation}.
  A \textit{function} is an operation that preserves the corresponding equalities.
  See~\cite{Pe21} for a more detailed explanation.}
$\chi_{\B A} : A^1 \cup A^0 \sto \D 2$, defined by
\[ \chi_{\B A}(x) := \left\{ \begin{array}{ll}
                 1   &\mbox{, $x \in A^1$}\\
                 0             &\mbox{, $x \in A^0$.}
                 \end{array}
          \right. \]
We call $\B A$ total, if $\dom(\B A) := A^1 \cup A^0 =_{\C P(X)} X$,
Let $\B A \subseteq \B B : \TOT A^1 \subseteq B^1 \ \& \ B^0 \subseteq A^0$,
and the totality of complemented subsets
$\C P^{\Disj}(X)$\index{$\C P^{\Disj}(X)$}\index{$\B A \subseteq \B B$} of $X$ is
equipped with the equality
$\B A =_{\C P^{\mathsmaller{\Disj}} (X)} \B B : \TOT \B A \subseteq \B B \ \& \ \B B \subseteq \B A$.
\end{defi}

Clearly, the complemented powerset $\C P^{\Disj}(X)$ of $X$ is a
proper class. If $f_1 \colon A^1 \subseteq B^1$ and $f_0 \colon B^0\subseteq A^0$,
then $f_1, f_0$ are strongly extensional functions.
E.g., if $f_1(a_1) \neq_{B^1} f_1(a_1{'})$, for some $a_1, a_1{'} \in A^1$,
then from the definition of the canonical inequality
$\neq_{B^1}$ we get $i_{B^1}\big(f_1(a_1)\big) \neq_{X} i_{B^1}\big(f_1(a_1{'})\big)$.
By the extensionality of $\neq_X$ we
get $i_{A^1}(a_1) \neq_X i_{A^1}(a_1{'}) :\TOT a_1 \neq_{A^1} a_1{'}$.

\begin{exa}\label{ex: detachable}
If $(X, =_X)$ is a set, let the following inequality on $X$:
\[x \neq_{(X, \mathsmaller{\D F(X, \D 2)})} x{'} : \TOT \exists_{f \in \D F(X, \D 2)}\big(f(x)
=_{\mathsmaller{\D 2}} 1 \ \& \ f(x{'}) =_{\mathsmaller{\D 2}} 0 \big)
\]
If $f \in \D F(X, \D 2)$, the following extensional subsets of $X$
\[ \delta_0^1(f) := \{x \in X \mid f(x) =_{\mathsmaller{\D 2}} 1\}, \]
\[ \delta_0^0(f) := \{x \in X \mid f(x) =_{\mathsmaller{\D 2}} 0\}, \]
are called detachable\index{detachable subset}, or free\index{free subset} subsets of $X$. Clearly,
$\B \delta (f) := \big(\delta_0^1(f), \delta_0^0(f)\big)$ is a complemented subset of $X$ with respect
to the inequality $\neq_X^{\mathsmaller{\D F(X, \D 2)}}$.
The characteristic function $\chi_{\B \delta(f)}$ of $\B \delta(f)$ is
(definitionally equal to) $f$ (recall that
$f(x) =_{\mathsmaller{\D 2}} 1 :\TOT f(x) := 1$), and $\delta_0^1(f) \cup \delta_0^0(f) = X$.
\end{exa}

\begin{rem}\label{rem: compl2}
 If $\B A \in \C P^{\Disj}(X)$, then $\B \chi_{\B A} := (A^1 \cup A^0, i_{A^1 \cup A^0}, \chi_{\B A}) \in
\F F^{\se}(X,\D 2)$.
\end{rem}

\begin{proof}
Let $z, w \in A^1 \cup A^0$
with $\chi_{\B A}(z) \neq_{\D 2} \chi_{\B A}(w)$. If for example
$\chi_{\B A}(z) := 1$
and $\chi_{\B A}(w) := 0$, then $z \in A^1$, $w \in A^0$. As $A^1 \Disj A^0$, we
get $i_{A^1}(z) \neq_X i_{A^0}(w) \TOT: z \neq_{A^1 \cup A^0} w$.
\end{proof}

\begin{defi}\label{def: foperationscomplemented2}
If $\B A, \B B \in \C P^{\Disj}(X)$, let the following operations\footnote{
  $\BMT$ and $\BCMT$ involve different operations on complemented subsets.
  We only describe the algebra of complemented subsets given in $\BCMT$,
  which has a more ``linear'' behavior (see also~\cite{Sh21}).
  For total complemented subsets the operations given in $\BMT$ and $\BCMT$ coincide.
  See~\cite{PW22} for an in-depth comparison of the two algebras of complemented subsets.}
on them:
$$\B A \vee \B B := \big([A^1 \cap B^1] \cup [A^1 \cap B^0] \cup [A^0 \cap B^1], \ A^0 \cap B^0\big),$$
$$\B A \wedge \B B := \big(A^1 \cap B^1, \ [A^1 \cap B^0] \cup [A^0 \cap B^1] \cup [A^0 \cap B^0]\big),$$
$$- \B A := (A^0, A^1),$$
$$\B A - \B B := \B A \wedge (- \B B).$$
\end{defi}

\begin{prop}\label{prp: second2}
$\big(\C P^{\Disj}(X), \wedge, \vee, -)$
satisfies all properties of a distributive lattice except\footnote{In~\cite{BB85}, p.~74,
  it is mentioned that complemented subsets satisfy
  ``all the usual finite algebraic laws that do not involve the operation of set complementation''.
  In~\cite{CP02}, p.~695, it is noticed though, that
  the absorption equalities are not satisfied.}
for the absorption
equalities $(\B A \wedge \B B) \vee \B A = \B A$ and $(\B A \vee \B B)
\wedge \B A = \B A$. Moreover, $- (- \B A) = \B A$, and $-(\B A \vee
\B B) = (-\B A) \wedge (-\B B)$, for every $\B A, \B B \in \C
P^{\Disj}(X)$.
\end{prop}

The classical bijection between $\C P(X)$ and $2^X$ is translated constructively as the existence of
``bijective'', proper class-assignment routines between the proper classes $\C P^{\Disj}(X) $ and
$\F F^{\se}(X, \D 2)$. The proof of the following fact is found in~\cite{PW22}, and it is the only place in
this paper that we refer to assignment routines defined on proper classes.

\begin{prop}\label{prp: compl3}
Consider the proper class-assignment routines
$$\chi \colon \C P^{\Disj}(X) \sto \F F^{\se}(X, \D 2) \ \ \& \ \
\delta \colon \F F^{\se}(X, \D 2) \sto \C P^{\Disj}(X),$$
$$\B A \mapsto \B \chi(\B A) \ \ \ \B f_A \mapsto \delta(\B f_A) := \big(\delta^1(f_A), \delta^0(f_A)\big),$$
\[ \delta^1(f_A) := \big\{a \in A \mid f_A(a) =_{\mathsmaller{\D 2}} 1 \big\} =: [f_A =_{\D 2} 1], \]
\[ \delta^0(f_A) := \big\{a \in A \mid f_A(a) =_{\mathsmaller{\D 2}} 0 \big\} =: [f_A =_{\D 2} 0]. \]
Then $\chi, \delta$ are well-defined, proper class-functions, which are
inverse to each other. Moreover, $\delta(\sim \B f) =_{\C P^{\Disj}(X)} - \delta(\B f)$ and
$\chi_{- \B A} =_{\F F(X, \D 2)} \sim
\chi_{\B A}$, where $\B f \in \F F^{\se}(X, \D 2)$ and $\B A \in \C P^{\Disj}(X)$.
\end{prop}

\begin{prop}\label{prp: second1}
Let $\B A, \B B \in \C P^{\Disj}(X)$ and $\B f, \B g \in \F F^{\se}(X, \D 2)$.
\begin{enumerate}[label=(\roman*)]
\item $\chi_{\B A \vee \B B} =_{\F F(X, \D 2)} \chi_{\B A} \vee \chi_{\B B}$,
  $\chi_{\B A \wedge \B B} =_{\F F(X, \D 2)} \chi_{\B A} \wedge \chi_{\B B}$.
\item $\chi_{- \B A} =_{\F F(X, \D 2)} 1 - \chi_{\B A}$
  and $\chi_{\B A - \B B} =_{\F F(X, \D 2)} \chi_{\B A}(1 - \chi_{\B B})$.
\item $ \delta ( \B f_A) \vee \delta (\B f_B) =_{\C P^{\Disj}(X)} \delta ( \B f_A \vee  \B f_B)$ and $\delta ( \B f_A) \wedge  \delta ( \B f_B) =_{\C P^{\Disj}(X)} \delta ( \B f_A \wedge  \B f_B)$.
\item $\delta ( \B f_A \sim \B f_B) =_{\C P^{\Disj}(X)} \delta ( \B f_A) - \delta(\B f_B)$.
\end{enumerate}
\end{prop}

\section{Families of subsets}
\label{sec: famsubsets}

In this section we present the basic notions and facts on set-indexed
families of subsets that are going to be used in the rest of the
paper.  Roughly speaking, a family of subsets of a set $X$ indexed by
some set $I$ is an assignment routine $\lambda_0 : I \sto \C P(X)$
that behaves like a function i.e., if $i =_I j$, then $\lambda_0(i)=_{\C P(X)} \lambda_0 (j)$.
The following definition is a formulation
of this rough description that reveals the witnesses of the equality
$\lambda_0(i) =_{\C P(X)} \lambda_0 (j)$. This is done ``internally'',
through the embeddings of the subsets into $X$. The equality
$\lambda_0(i) =_{\D V_0} \lambda_0 (j)$, which is defined
``externally'' through the transport maps (see~\cite{Pe21}, Definition
3.1), follows, and a family of subsets is also a family of sets. We
start by introducing some notation. For details we refer
to~\cite{Pe21}.

\begin{defi}\label{def: depoperation}
  Let $I$ be as set and $\lambda_0 : I \sto \D V_0$. A \emph{dependent operation} over $\lambda_0$
  \[
  \Phi : \bigcurlywedge_{i \in I}\lambda_0(i)
  \]
  assigns to each $i\in I$ an element $\Phi(i):= \Phi_i\in\lambda_0(i)$. We denote by $\D A(I, \lambda_0)$
  the totality of dependent operations over $\lambda_0$ equipped with the equality
  \[
  \Phi =_{\D A(I, \lambda_0)} \Psi :\Leftrightarrow \forall_{i\in I} \big(\Phi_i =_{\lambda_0(i)}\Psi_i\big).
  \]
\end{defi}

\begin{defi}\label{def: famofsubsets}
Let $X$ and $I$ be sets and let $D(I) := \{(i, i{'}) \in I \times I \mid i =_I i{'}\}$ be the diagonal of $I$.
A \textit{family of subsets}\index{family of subsets} of $X$ indexed by $I$,
is a triplet
$\Lambda(X) := (\lambda_0, \C E, \lambda_1)$, where
\index{$\Lambda(X)$}\index{$\C E$}
$\lambda_0 : I \sto \D V_0$,
\[ \C E : \bigcurlywedge_{i \in I}\D F\big(\lambda_0(i), X\big), \ \ \ \ \C E(i) := \C E_i; \ \ \ \ i \in I, \]
\[ \lambda_1 : \bigcurlywedge_{(i, j) \in D(I)}\D F\big(\lambda_0(i), \lambda_0(j)\big), \ \ \ \
\lambda_1(i, j) := \lambda_{ij}; \ \ \ \ (i, j) \in D(I), \]
such that the following conditions hold:
\begin{enumerate}[label=(\roman*)]
\item For every $i \in I$, the function $\C E_i : \lambda_0(i) \to X$ is an embedding.
\item For every $i \in I$, we have that $\lambda_{ii} =_{\D F(\lambda_0(i), \lambda_0(i))} \id_{\lambda_0(i)}$.
\item For every $(i, j) \in D(I)$ we have that
$\C E_i =_{\D F(\lambda_0(i), X)} \C E_j \circ \lambda_{ij}$.
\end{enumerate}
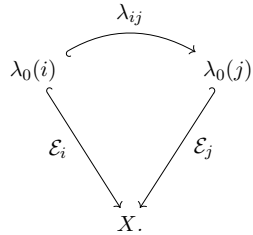
\begin{center}
\resizebox{3.5cm}{!}{%
\begin{tikzpicture}

\node (E) at (0,0) {$\lambda_0(i)$};
\node[right=of E] (B) {};
\node[right=of B] (F) {$\lambda_0(j)$};
\node[below=of B] (C) {};
\node[below=of C] (A) {$X$.};

\draw[left hook->,bend left] (E) to node [midway,above] {$\lambda_{ij}$} (F);
\draw[right hook->] (E)--(A) node [midway,left] {$\C E_{i} \ $};
\draw[left hook->] (F)--(A) node [midway,right] {$ \ \C E_j$};

\end{tikzpicture}
}
\end{center}
We call a pair $A_i := (\lambda_0(i), \C E_i)$ an element of $\Lambda(X)$.
If $(A, i_A) \in \C P(X)$, the \textit{constant} $I$-\textit{family of subsets}\index{constant family of subsets}
$A$\index{$C^A_X$}
is the pair
$C^{A}(X) := (\lambda_0^{A}, \C E^{A}, \lambda_1^A)$, where $\lambda_0 (i) := A$, $\C E_i^{A} := i_A$, and
$\lambda_1 (i, j) := \id_A$, for every $i \in I$ and $(i, j) \in D(I)$, respectively.
If $(A, i_A), (B, i_B) \subseteq X$, the
triplet \index{$\Lambda^{\D 2}(X)$}$\Lambda^{\D 2}(X) := (\lambda_0^{\D 2}, \C E, \lambda_1^{\D 2})$, where
$\lambda_0^{\D 2}(0):=A$ and $\lambda_0^{\D 2}(1):=B$,
$\C E_0 := i_A$ and $\C E_1 := i_B$,
$\lambda_1^{\D 2}(0,0):=\id_A$ and $\lambda_1^{\D 2}(1,1):=\id_B$
is the $\D 2$-family of subsets $A$ and $B$\index{$\D 2$-family of subsets} of $X$.
If $\Fam(I,X)$\index{$\Fam(I,X)$} denotes
the totality of $I$-families of subsets of $X$,
its equality is defined as in~\cite{Pe21}, Definition 3.2.
\end{defi}

\begin{exa}\label{ex: detachablefam}
Let $\big(X, =_X, \neq_{(X, \mathsmaller{\D F(X, \D 2)})}\big)$ be the set with inequality from
Example~\ref{ex: detachable}. The family of subsets
$\Delta^1(X) := \big(\delta_0^1, \C E^{1}, \delta_1^1\big)$ over the index-set $\D F(X, 2)$ is defined by the following rules:
$$\delta_0^1 \colon \D F(X, \D 2) \sto \D V_0, \ \ \ \ f \mapsto \delta_0^1(f), \ \ \ f \in \D F(X, \D 2),$$
$$\C E^{1} \colon \bigcurlywedge_{f \in \D F(X, \D 2)}\D F(\delta_0^1(f), X), \ \ \
\ \C E^{1}_f \colon \delta_0^1(f) \eto X \ \ \ \ x \mapsto x; \ \ \ \ x \in \delta_0^1(f),$$
$$\delta_1^1 \colon \bigcurlywedge_{(f,g) \in D(\D F(X, \D 2))}\D F(\delta_0^1(f), \delta_0^1(g)), \ \ \ \
\delta_1^1(f,g) := \delta_{fg}^1 \colon \delta_0^1(f) \to \delta_0^1(g) \ \ \ \ x \mapsto x; \ \ \ \
x \in \delta_0^1(f).$$
If $\Delta^0(X) := \big(\delta_0^0, \C E^{0}, \delta_1^0\big)$, where\index{$\Delta^0(X)$}
$\delta_0^0 \colon \D F(X, \D 2) \sto \D V_0$ is defined by the rule $f \mapsto \delta_0^0(f)$,
for every $f \in \D F(X, \D 2)$, and the dependent operations $\C E^{0}, \delta_1^0$ are defined similarly
to $\C E^{1}$ and $\delta_1^1$,
then $\Delta^1(X), \Delta^0(X)$ are \textit{sets} of subsets of $X$ in the following sense.
\end{exa}

\begin{defi}\label{def: setofsubsets}
Let $\Lambda(X) := (\lambda_0, \C E, \lambda_1) \in \Fam(I, X)$. We say that $\Lambda(X)$ is a
\emph{set} of subsets of $X$ if
\[
\forall_{i,j\in I}\big(\lambda_0(i) =_{\C P(X)} \lambda_0(j) \To i=_I j\big).
\]
In this case we write $\Lambda(X)\in\Set(I,X)$.
We can always make $\Lambda(X)$ into a set of subsets
$\tilde{\Lambda}(X) := (\lambda_0, \C E, \tilde\lambda_1)$
indexed by the set $\lambda_0 I(X)$, where $\lambda_0 I(X)$
is the totality $I$ with a new equality given by
\[ i =_{\lambda_0 I(X)} j :\TOT \lambda_0(i) =_{\C P(X)}\lambda_0(j),\]
for every $i, j \in I$. The assignment routine $\lambda_0:I\sto\D V_0$ and the dependent
function $\C E : \bigcurlywedge_{i \in I}\D F\big(\lambda_0(i), X\big)$
are the same as in $\Lambda (X)$.
Using the \textit{dependent version of Myhill's axiom of unique choice}\footnote{According to it, if $(\mu_0, \mu_1)$ is an $I$-family of sets such that
  for every $i \in I$ there is a unique (up to equality) $x_i \in \mu_0(i)$,
  then there is a dependent
  assignment routine $\Phi \colon \bigcurlywedge_{i \in I} \mu_0(i)$.
  The non-dependent version of this axiom is generally
  accepted by the practitioners of $\BISH$ and it is included in
  Myhill's system $\CST$ in~\cite{My75}.
  If $i =_{\lambda_0 I(X)} j$, then by the equality
  $\lambda_0(i) =_{\C P(X)}\lambda_0(j)$
  there is a \emph{unique} function $\lambda_0(i)\to \lambda_0(j)$
  commuting with the embeddings
  $\C E_i$ and $\C E_j$.
  To avoid Myhill's axiom, we need to add to
  our data a dependent assignment
  routine $\Delta$ that corresponds to
  every element of the diagonal of $\lambda_0 I(X)$
  an element of the set of witnesses $(e, e{'})$ of
  the corresponding equalities.},
one can define the dependent function
$\tilde\lambda_1:\bigcurlywedge_{(i, j) \in D(\lambda_0 I(X))}\D F\big(\lambda_0(i), \lambda_0(j)\big)$.
\end{defi}

As we explained in~\cite{Pe21}, the totality $\Fam(I)$ of all
$I$-families of sets cannot be accepted as a set, as the constant
$I$-family with value $\Fam(I)$ would then be defined through a
totality in which it belongs to.  This does not work as an argument
against the set-character of $\Fam(I, X)$.  It is not clear how the
constant $I$-family $\Fam(I, X)$ can be seen as a family of subsets of
$X$. If $\nu_0(i) := \Fam(I, X)$, for every $i \in I$, we need to
define a modulus of embeddings $\C N_i \colon \Fam(I, X) \eto X$, for
every $i \in I$. From the given data one could define the assignment
routine $\C N_i$ by the rule $\C N_i\big(\Lambda(X)\big) := \C E_i(u_i)$,
if it is known that $u_i \in \lambda_0(i)$. Even in that
case, the assignment routine $\C N_i$ cannot be shown to satisfy the
expected properties. Clearly, if $\C N_i$ was defined by the rule
$\C N_i\big(\Lambda(X)\big) := x_0 \in X$, then it cannot be an
embedding. The set-character of the totality $\Fam(I, X)$ is related
to the definition of a pre-measure space (see also the discussion
after the definition of a Bishop-Cheng measure space in
section~\ref{sec: premeasures}).  Next we describe the Sigma- and the
Pi-set of a family of subsets.

\begin{defi}\label{def: interiorunion}
Let $\Lambda(X) := (\lambda_0, \C E, \lambda_1) \in \Fam(I, X)$.
The \textit{interior union}\index{interior union}, or the union\index{union of a family of subsets}
of $\Lambda(X)$ is the totality $\sum_{i \in I} \lambda_0(i)$, which we
denote in this case by $\bigcup_{i \in I} \lambda_0(i)$.
Let the non-dependent assignment routine
$e \colon \bigcup_{i \in I} \lambda_0(i) \sto X$ defined by
$(i, x) \mapsto \C E_i (x)$, for every $(i, x) \in \bigcup_{i \in I} \lambda_0(i)$, and let
\[ (i, x) =_{\mathsmaller{\bigcup_{i \in I} \lambda_0(i)}} (j, y) :\Leftrightarrow
e(i, x) =_X
e(j, y) :\Leftrightarrow \C E_i (x) =_X \C E_j (y). \]
If $\neq_X$ is an inequality on $X$, let
$(i, x) \neq_{\mathsmaller{\bigcup_{i \in I} \lambda_0(i)}} (j, y) :\Leftrightarrow
\C E_i (x) \neq_X \C E_j (y)$.
The family $\Lambda(X)$ is called a \textit{covering}\index{covering} of $X$, or
$\Lambda(X)$ covers $X$\index{$\Lambda(X)$ covers $X$},  if
$\bigcup_{i \in I} \lambda_0(i)  =_{\C P (X)} X$.
If $\neq_I$ is an inequality on $I$, we say that $\Lambda(X)$ is a
family of disjoint subsets of $X$\index{family of disjoint subsets of $X$} $($with respect to $\neq_I)$,  if
$\forall_{i, j \in I}\big(i \neq_I j \To \lambda_0(i) \Disj \lambda_0(j)\big)$,
where by Definition~\ref{def: apartsubsets}
$\lambda_0(i) \Disj \lambda_0(j) :\TOT \forall_{u \in \lambda_0(i)}\forall_{w \in
\lambda_0(j)}\big(\C E_i(u) \neq_X \C E_j(w)\big)$.
\end{defi}

Clearly, $=_{\mathsmaller{\bigcup_{i \in I} \lambda_0(i)}}$ is an equality on
$\bigcup_{i \in I} \lambda_0(i)$,
and the
operation $e$ is an embedding of $\bigcup_{i \in I} \lambda_0(i)$ into $X$, hence
$\big(\bigcup_{i \in I} \lambda_0(i), e\big) \subseteq X$. The inequality
$\neq_{\mathsmaller{\bigcup_{i \in I} \lambda_0(i)}}$ is the canonical inequality of the subset
$\bigcup_{i \in I} \lambda_0(i)$ of $X$.
Hence, if $(X, =_X, \neq_X)$ is discrete, then $\big(\bigcup_{i \in I} \lambda_0(i),
=_{\mathsmaller{\bigcup_{i \in I} \lambda_0(i)}}, \neq_{\mathsmaller{\bigcup_{i \in I} \lambda_0(i)}}\big)$ is discrete,
and if $\neq_X$ is tight, then $\neq_{\mathsmaller{\bigcup_{i \in I} \lambda_0(i)}}$ is tight.
As the following left diagram commutes, $\Lambda(X)$ covers $X$, if and only if the following right
diagram commutes i.e., if and only if $X \subseteq \bigcup_{i \in I} \lambda_0(i)$
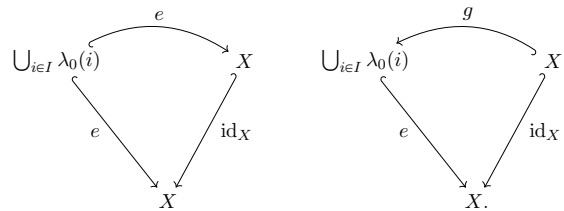
\begin{center}
\resizebox{7.5cm}{!}{%
\begin{tikzpicture}

\node (E) at (0,0) {$\bigcup_{i \in I} \lambda_0(i)$};
\node[right=of E] (B) {};
\node[right=of B] (F) {$X$};
\node[below=of B] (C) {};
\node[below=of C] (A) {$X$};

\node[right=of F] (G) {$\bigcup_{i \in I} \lambda_0(i)$};
\node[right=of G] (H) {};
\node[right=of H] (J) {$X$};
\node[below=of H] (K) {};
\node[below=of K] (L) {$X$.};

\draw[left hook->,bend left] (E) to node [midway,above] {$e$} (F);
\draw[right hook->] (E)--(A) node [midway,left] {$e \ $};
\draw[left hook->] (F)--(A) node [midway,right] {$ \ \id_X$};
\draw[right hook->,bend right] (J) to node [midway,above] {$g$} (G);
\draw[right hook->] (G)--(L) node [midway,left] {$e \ $};
\draw[left hook->] (J)--(L) node [midway,right] {$ \ \id_X$};

\end{tikzpicture}
}
\end{center}

If $(i, x) =_{\mathsmaller{\bigcup_{i \in I} \lambda_0(i)}} (j, y)$, it is not
necessary that $i =_I j$, hence it is not necessary that
$(i, x) =_{\mathsmaller{\sum_{i \in I} \lambda_0(i)}} (j, y)$ (as we show in the next proposition, the converse
implication holds). Consequently, the first projection operation $\pr_1^{\Lambda(X)}
:= \pr_1^{\Lambda}$\index{$\pr_1^{\Lambda(X)}$}, where $\Lambda$ is the $I$-family of sets induced by $\Lambda(X)$,
is not necessarily a function!
The second projection map on $\Lambda(X)$ is defined by
$\pr_2^{\Lambda(X)} := \pr_2^{\Lambda}$\index{$\pr_2^{\Lambda(X)}$}.
Notice that $\neq_{\mathsmaller{\bigcup_{i \in I} \lambda_0(i)}}$ is an
inequality on $\bigcup_{i \in I} \lambda_0(i)$, without supposing neither an inequality on $I$, nor an inequality on the sets $\lambda_0(i)$'s.
The following remarks are straightforward to show.

\begin{rem}\label{remp: interiorunion1}
Let $\Lambda(X) := (\lambda_0, \C E, \lambda_1) \in \Fam(I, X)$.
\begin{enumerate}[label=(\roman*)]
\item If $(i, x) =_{\mathsmaller{\sum_{i \in I} \lambda_0(i)}} (j, y)$, then
  $(i, x) =_{\mathsmaller{\bigcup_{i \in I} \lambda_0(i)}} (j, y)$.
\item If $e \colon \sum_{i \in I} \lambda_0(i) \sto X$ is an embedding, then
$\big(\sum_{i \in I} \lambda_0(i), e\big) =_{\C P(X)} \big(\bigcup_{i \in I} \lambda_0(i),e\big)$.
\item If $\neq_I$ is tight and $\Lambda(X)$ is a family of disjoint subsets
with respect to $\neq_I$, then  $e \colon \sum_{i \in I} \lambda_0(i) \eto X$.
\end{enumerate}
\end{rem}

\begin{rem}\label{rem: unionA}
If $i_0 \in I$, $(A, i_A) \subseteq X$, and $C^A(X) := (\lambda_0^A, \C E^{A}, \lambda_1^A)
\in \Fam(I, X)$ is the constant family $A$ of subsets of $X$, then
\[ \bigcup_{i \in I}A := \bigcup_{i \in I}\lambda_0^A(i) =_{\C P(X)} A. \]

\end{rem}

\begin{rem}\label{rem: interiorunion2}
If $\Lambda^{\D 2}(X)$ is the $\D 2$-family of subsets $A, B$ of $X$,
$\bigcup_{i \in \D 2} \lambda_0^{\D 2}(i) =_{\C P(X)} A \cup B.$
\end{rem}

\begin{defi}\label{def: intfamilyofsubsets}
Let $\Lambda(X) := (\lambda_0, \C E, \lambda_1) \in \Fam(I, X)$, and $i_0 \in I$.
The intersection $\bigcap_{i \in I} \lambda_0 (i)$\index{intersection of a family of
subsets}\index{$\bigcap_{i \in I} \lambda_0 (i)$} of $\Lambda(X)$ is the
totality defined by
\[ \Phi \in \bigcap_{i \in I} \lambda_0 (i) : \TOT \Phi \in \D A(I, \lambda_0)
 \ \& \ \forall_{i, j \in I}\big(\C E_i(\Phi_i) =_X \C E_{j} (\Phi_{j})\big). \]
Let $e \colon \bigcap_{i \in I} \lambda_0 (i) \sto X$ be defined by
$e(\Phi) := \C E_{i_0}\big(\Phi_{i_0}\big)$, for
every $\Phi \in \bigcap_{i \in I} \lambda_0 (i)$, and
\[ \Phi =_{\mathsmaller{\bigcap_{i \in I} \lambda_0 (i)}} \Theta : \TOT e(\Phi)
=_X e(\Theta) :\TOT
\C E_{i_0}\big(\Phi_{i_0}\big) =_X \C E_{i_0}\big(\Theta_{i_0}\big),
\]
If $\neq_X$ is a given inequality on $X$, let
$\Phi \neq_{\mathsmaller{\bigcap_{i \in I} \lambda_0 (i)}} \Theta :\TOT
 \C E_{i_0}\big(\Phi_{i_0}\big) \neq_X \C E_{i_0}\big(\Theta_{i_0}\big)$.
\end{defi}

The following remarks are straightforward to show.

\begin{rem}\label{rem: intersection1}
  Let $\Lambda(X) := (\lambda_0, \C E, \lambda_1) \in \Fam(I, X)$.
  \begin{enumerate}[label=(\roman*)]
  \item $\Phi =_{\bigcap_{i \in I} \lambda_0 (i)} \Theta \TOT \Phi =_{\D A(I, \lambda_0)} \Theta$.
  \item If $\Phi \in \bigcap_{i \in I} \lambda_0 (i)$, then $\Phi \in \prod_{i \in I}\lambda_0(i)$.
  \item If $(X, =_X, \neq_X)$ is discrete, the set $\big(\bigcap_{i \in I} \lambda_0(i),
=_{\mathsmaller{\bigcap_{i \in I} \lambda_0(i)}}, \neq_{\mathsmaller{\bigcap_{i \in I} \lambda_0(i)}}\big)$ is discrete.
  \end{enumerate}
\end{rem}

\begin{rem}\label{rem: intA}
Let $i_0 \in I$, $(A, i_A) \subseteq X$, and $C^A(X) := (\lambda_0^A, \C E^{A}, \lambda_1^A) \in \Fam(I, X)$
the constant family $A$ of subsets of $X$. Then
\[ \bigcap_{i \in I}A := \bigcap_{i \in I}\lambda_0^A(i) =_{\C P(X)} A. \]

\end{rem}

\begin{rem}\label{rem: intersection2}
If $\Lambda^{\D 2}(X)$ is the $\D 2$-family of subsets $A, B$ of $X$,
$\bigcap_{i \in \D 2} \lambda_0^{\D 2}(i) =_{\C P(X)} A \cap B.$
\end{rem}

\section{Families of partial functions and families of complemented subsets}
\label{sec: fampartial}

Next we define within $\BST$ the notions of a family of partial
functions and of a family of complemented subsets indexed by some set
$I$. These function-theoretic concepts will be used in $\PBCMT$
instead of the abstract sets of partial functions and of subsets,
respectively, that are considered in $\BCMT$.

\begin{defi}\label{def: famofpartial}
Let $X, Y$ and $I$ be sets. A \textit{family of partial functions}\index{family of partial functions}
from $X$ to $Y$ indexed by $I$,
or an $I$-\textit{family of partial functions from $X$ to $Y$}\index{$I$-family of partial functions}, is a triplet
$\Lambda(X,Y) := (\lambda_0, \C E, \lambda_1, \F f)$, where\index{$\Lambda(X,Y)$}\index{$\F f^Y$}
$\Lambda(X) := (\lambda_0, \C E, \lambda_1) \in \Fam(I, X)$ and
$\F f : \bigcurlywedge_{i \in I}\D F\big(\lambda_0(i), Y\big)$ with
$\F f(i) := \F f_i$, for every $i \in I$,
such that, for every $(i, j) \in D(I)$, the following
diagrams commute
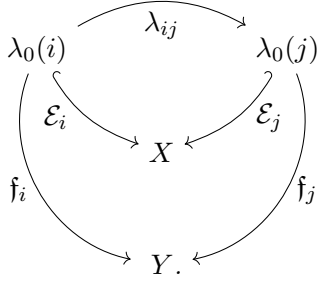
\begin{center}
\resizebox{4.5cm}{!}{%
\begin{tikzpicture}

\node (E) at (0,0) {$\lambda_0(i)$};
\node[right=of E] (B) {};
\node[right=of B] (F) {$\lambda_0(j)$};
\node[below=of B] (A) {$X$};
\node[below=of A] (C) {$ \ Y$.};

\draw[->,bend left] (E) to node [midway,below] {$\lambda_{ij}$} (F);
\draw[->,bend right=50] (E) to node [midway,left] {$\F f_i$} (C);
\draw[->,bend left=50] (F) to node [midway,right] {$\F f_j$} (C);
\draw[right hook->,bend right=20] (E) to node [midway,left] {$ \ \C E_i \ $} (A);
\draw[left hook->,bend left=20] (F) to node [midway,right] {$\ \C E_j \ $} (A);

\end{tikzpicture}
}
\end{center}
If $i \in I$, we call the partial function
$\B f_i := (\lambda_0(i), \C E_i, \F f_i) \in \F F(X, Y)$
an element of $\Lambda(X,Y)$.

The equality on the totality $\Fam(I,X,Y)$\index{$\Fam(I,X,Y)$}
of $I$-families of partial functions from $X$ to $Y$ can be defined in an obvious way, analogously
to the equality on $\Fam(I,X)$ given in~\cite[Def. 3.2]{Pe21}.
\end{defi}

Clearly, if $\Lambda(X,Y) \in \Fam(I,X,Y)$ and $(i,j) \in D(I)$, then $(\lambda_{ij}, \lambda_{ji}) \colon
\B f_i =_{\mathsmaller{\F F(X,Y)}} \B f_j$.

\begin{defi}\label{def: setofpartialfunctions}
Let $\Lambda(X,Y) := (\lambda_0, \C E, \lambda_1,\F f) \in \Fam(I, X,Y)$. We say that $\Lambda(X,Y)$
is a \emph{set} of partial functions from $X$ to $Y$ if
\[
\forall_{i,j\in I}\big(\B f_i =_{\F F(X,Y)} \B f_j \To i=_I j\big).
\]
In this case we write $\Lambda(X,Y)\in\Set(I,X,Y)$\index{$\Set(I,X,Y)$} and even
$\Lambda(X,Y)\in\Set^\se(I,X,Y)$\index{$\Set^\se(I,X,Y)$} if $\Lambda(X,Y)$ is a family
of strongly extensional partial functions.

As described in Definition~\ref{def: setofsubsets},
we can make $\Lambda(X,Y)$ into a $\lambda_0I(X,Y)$-set of partial functions
$\tilde\Lambda(X,Y) := (\lambda_0, \C E, \tilde\lambda_1,\F f)$.
As in the case for subsets, $\lambda_0I(X,Y)$ is the totality $I$ equipped with the equality
\[ i =_{\lambda_0I (X, Y)} j :\TOT \B f_i =_{\F F(X,Y)} \B f_j,\]
or if $\Lambda(X,Y)$ is a family of strongly extensional partial functions,
\[ i =_{\lambda_0I (X, Y)} j :\TOT \B f_i =_{\F F^{\se}(X,Y)} \B f_j.\]
The only component changing is the new $\tilde\lambda_1$, which is defined using
dependent unique choice.
\end{defi}

\begin{defi}\label{def: famofcompl}
Let the sets $(X, =_X, \neq_X)$ and $(I, =_I)$.
A \textit{family of complemented subsets}\index{family of complemented subsets}
of $X$ indexed by $I$, or an $I$-\textit{family of complemented subsets of $X$}\index{$I$-family of complemented subsets},
is a structure\index{$\B \Lambda(X)$}
$\B \Lambda(X) := \big(\lambda_0^1, \C E^{1}, \lambda_1^1, \lambda_0^0, \C E^{0}, \lambda_1^0)$, such that
$\Lambda^1(X) := \big(\lambda_0^1, \C E^{1}, \lambda_1^1\big) \in \Fam(I, X)$ and
$\Lambda^0(X) := \big(\lambda_0^0, \C E^{0}, \lambda_1^0\big) \in \Fam(I, X)$
i.e., for every $(i, j) \in D(I)$, the following diagrams commute
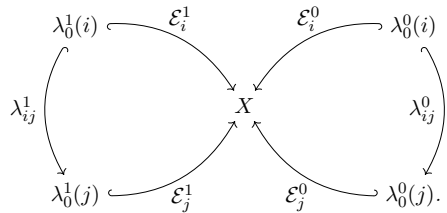
\begin{center}
\resizebox{6cm}{!}{%
\begin{tikzpicture}

\node (E) at (0,0) {$X$};
\node[above=of E] (F) {};
\node[right=of F] (P) {};
\node[right=of P] (G) {$\lambda_0^0(i)$};
\node[below=of G] (S) {};
\node[below=of S] (H) {$\lambda_0^0(j)$.};
\node[left=of F] (X) {};
\node[left=of X] (J) {$\lambda_0^1(i)$};
\node[below=of J] (T) {};
\node[below=of T] (K) {$\lambda_0^1(j)$};

\draw[right hook->,bend right] (K) to node [midway,below] {$\C E_j^{1}$} (E);
\draw[left hook->,bend left] (J) to node [midway,above] {$\C E_i^{1}$} (E);
\draw[right hook->,bend right] (G) to node [midway,above] {$\C E_i^{0}$} (E);
\draw[left hook->,bend left] (H) to node [midway,below] {$\C E_j^{0} \ \ $} (E);
\draw[right hook->,bend right] (J) to node [midway,left] {$\lambda_{ij}^1$} (K);
\draw[left hook->,bend left] (G) to node [midway,left] {$\lambda_{ij}^0$} (H);

\end{tikzpicture}
}
\end{center}
Moreover, for every $i \in I$ the element
$\B \lambda_0(i) := \big(\lambda_0^1(i), \lambda_0^0(i)\big)$ of $\B \Lambda(X)$ is in $\C P^{\Disj}(X)$.
Again, the equality on $\Fam(I, \B X)$, the totality of $I$-families of complemented subsets
of $X$, is defined in an obvious way, analogously to \cite{Pe21}, Definition 3.2.
\end{defi}

As in the case of $\Fam(I, X)$, we assume the totality $\Fam(I, \B X)$ to be a set.
The operations
$\wedge$ and $\vee$ between complemented subsets in Definition~\ref{def: foperationscomplemented2} are
extended to families of complemented subsets.
We write $\Set(I, \B X)$\index{$\Set(I, \B X)$}
for the totality of $I$-\emph{sets} of complemented subsets of $X$,
which are defined completely analogously to sets of subsets or partial functions.

\section{Pre-measure spaces}
\label{sec: premeasures}

In this section we introduce the notion of a pre-measure space as a
predicative counterpart to the notion of Bishop-Cheng measure space in
$\BCMT$. The pre-measure space of complemented detachable subsets of a
set $X$ with the Dirac-measure concentrated at a single point is
described. The notion of a Bishop-Cheng measure space is defined
in~\cite{BB85}, p.~282, and
appeared first\footnote{In~\cite{BC72}, p.~55,
  condition $(\BCMS_1)$ appears in the equivalent form: if $\B B \in M$ with
  $B^1 \subseteq A^1$ and $B^0 \subseteq A^0$, then $\B A \in M$.}
in~\cite{BC72} p.~47.

\begin{defiC}[\textbf{Bishop-Cheng measure space}]
A (Bishop-Cheng) measure space is a triplet $(X, M, \mu)$ consisting of an inhabited
set with inequality $(X, =_X, \neq_X)$, a set $M$ of complemented sets in
$X$, and a mapping $\mu$ of $M$ into $\Real^{0+}$, such that the following properties
hold:
\begin{enumerate}[label=$(\BCMS_{\arabic*})$]
\item If $\B A$ and $\B B$ belong to $M$, then so do $\B A \vee \B B$ and $\B  A \wedge \B  B$, and $\mu(\B A) + \mu(\B B) = \mu(\B A \vee \B B) + \mu(\B A \wedge \B B)$.
\item If $\B A$ and $\B A \wedge \B B$ belong to $M$, then so does $\B A - \B B$, and $\mu(\B A) = \mu(\B A \wedge \B B) + \mu(\B A - \B B)$.
\item There exists $\B A$ in $M$ such that $\mu(\B A) >0$.
\item If $(\B A_n)$ is a sequence of elements of $M$ such that
  $\lim_{k \to \infty}\mu\big(\bigwedge_{n = 1}^k \B A_n\big)$ exists and is positive, then $\bigcap_{n}\B A_n^1$ is inhabited.
\end{enumerate}
\end{defiC}
\noindent
The elements of $M$ are the \textit{integrable sets} of the
measure space $(X, M, \mu)$, and for each $\B A$ in $M$ the non negative number $\mu(\B A)$ is the
\textit{measure} of $\B A$. In the above definition there is no indication how the set $M$ of complemented sets
of $X$ is constructed, and $(\BCMS_2)$ requires quantification over the universe $\D V_0$, as $\B B$ is an arbitrary
complemented subsets of $X$. In~\cite{Bi67}, p.~183, Bishop used two families of complemented subsets,
in order to avoid such a quantification in his definition of a measure space within $\BMT$. One set-indexed
family which $\B A$ and $\B A \wedge \B B$ belong to, and one which $\B B$ belongs to.
In Definition~\ref{def: premeasurespace} we predicatively reformulate the Bishop-Cheng definition of
a measure space.
Especially for condition $(\BCMS_2)$ we provide two alternatives.
In the first one, condition $(\PMS_2)$ in Definition~\ref{def: premeasurespace}, we use the fact that within $\BST$ the totality
$\Fam(\D 1, X)$ of $\D 1$-families of complemented subsets, where $\D 1 := \{0\}$, is assumed to be a set\footnote{Notice that
in order to define an $\D 1$-family of complemented subsets we need \textit{first}
to construct a complemented subset $(A^1, A^0)$
of $X$, and \textit{afterwards} to define $\lambda_0^0(0) := A^1$ and $\lambda_0^0(0) := A^0$.}, hence
quantification over it is allowed.
In the second alternative, the weaker condition $(\PMS_2^*)$ in Definition~\ref{def: premeasurespace},
only quantification over the index-set is required.
If $\B \Lambda(X)$ is an $I$-family of complemented subsets of $X$,
we can define an equality on the index set $I$, such that
the converse implication $\B \lambda_0(i) =_{\C P^{\Disj}(X)} \B \lambda_0(j) \To i = j$ also holds.
The family $\B \Lambda(X)$ is then called (as in the case of a family of subsets in~\cite{Bi67}, p.~65)
a \textit{set} of complemented subsets. Consequently, functions
on the index-set $I$ are extended to functions on the complemented subsets $\B \lambda_0(i)$'s.

One could predicatively reformulate the definition of a Bishop-Cheng
measure space within $\BST$. We proceed instead directly to define the
notion of a pre-measure space, giving an explicit formulation of
Bishop's suggestion, expressed in~\cite{Bi70}, p.~67, with respect to
Definition~\ref{def: premeasurespace}. The main idea is to define
operations on $I$ that correspond to the operations on complemented
subsets, and reformulate accordingly the clauses for the measure
$\mu$.

\begin{defi}[Pre-measure space within $\BST$]\label{def: premeasurespace}
Let $(X, =_X, \neq_X)$ be an inhabited set, and let $(I, =_I)$ be
equipped with operations $\vee \colon I \times I \sto I, \wedge \colon I \times I \sto I$ $($for
simplicity we use the same
symbols with the ones for the operations on complemented subsets$)$,
and $\sim \colon I\times I \sto I$.
If $i, j \in I$, let $i \leq j : \TOT i \wedge j = i$.
Let $\B \Lambda(X) := (\lambda_0^1, \C E^{1}, \lambda_1^1 ; \lambda_0^0, \C E^{0}, \lambda_1^0)
\in \Set(I, \B X)$, and $\mu \colon I \to [0, + \infty)$ such that the following conditions hold:
\[ (\PMS_1) \  \  \forall_{i, j \in I}\bigg(\B \lambda_0 (i) \vee \B \lambda_0 (j) = \B \lambda_0 (i \vee j)
\ \& \B \lambda_0 (i) \wedge \B \lambda_0 (j) = \B \lambda_0 (i \wedge j) \ \& \ \ \ \ \ \ \ \ \ \ \ \ \ \ \ \ \ \ \ \ \ \ \ \ \ \ \ \]
\[ \ \ \ \ \ \ \ \B \lambda_0 (i) - \B \lambda_0 (j) =
\B \lambda_0 (i \sim j) \ \& \ \mu(i) + \mu(j) = \mu(i \vee j) + \mu(i \wedge j)\bigg).\]
\[ (\PMS_2) \ \ \   \ \ \ \ \ \ \ \ \ \ \ \ \ \ \ \ \   \ \ \ \ \ \
\forall_{i \in I}\forall_{\B A(X) \in \Fam(\D 1, \B X)}\bigg[\exists_{k \in I}\bigg(\B \lambda_0 (i) \
\wedge \ \B \alpha_0 (0) = \B \lambda_0 (k)\bigg) \ \To \ \ \ \ \ \ \ \ \ \ \ \ \ \ \ \
\ \ \ \  \ \ \ \ \ \ \ \ \ \ \ \  \]
\[ \ \ \ \ \ \ \ \ \ \ \ \ \ \B \lambda_0 (i) - \B \alpha_0(0) = \B \lambda_0 (i \sim k) \ \& \ \mu(i) = \mu(k) + \mu(i \sim k)\bigg]. \ \ \ \ \ \ \ \ \ \ \]
\[ (\PMS_3) \ \ \ \ \ \ \ \ \ \ \  \ \ \ \ \ \  \ \ \  \ \ \ \ \ \ \ \ \ \ \ \ \ \  \ \ \ \ \ \ \ \ \ \ \ \ \ \ \ \ \
\exists_{i \in I}\big(\mu (i)\big) > 0. \ \ \ \ \ \ \ \  \ \ \ \ \ \ \ \ \ \ \  \ \ \ \ \ \ \ \  \
\ \ \ \ \ \ \ \ \ \ \ \ \ \ \ \ \ \ \ \ \ \ \ \ \ \  \]
\[ (\PMS_4) \ \ \ \ \ \ \ \ \  \ \ \ \ \ \ \ \ \ \ \ \ \  \forall_{\alpha \in \D F(\Nat, I)}\bigg[\exists
\lim_{\mathsmaller{m \to + \infty}} \mu \bigg(\bigwedge_{n = 1}^m \alpha(n)\bigg) \  \& \
\lim_{\mathsmaller{m \to + \infty}} \mu \bigg(\bigwedge_{n = 1}^m \alpha(n)\bigg) > 0 \To \ \ \ \ \
\ \ \ \  \ \ \ \ \ \ \ \ \ \]
\[ \ \ \ \ \ \  \ \ \ \ \  \exists_{x \in X}\bigg(x \in \bigcap_{n \in \Nat}\lambda_0^1 (\alpha (n))\bigg)\bigg].
\ \ \ \ \ \  \ \]
We call the triplet $\C M (\B \Lambda(X)) := (X, I, \mu)$
a \textit{pre-measure space}\index{pre-measure space}, the function $\mu$ a \textit{pre-measure}\index{pre-measure},
and the index-set $I$ the set of integrable, or measurable \index{set of integrable indices}indices of
$\C M (\B \Lambda(X))$.

Alternatively to $(\PMS_2)$, one could use the following clause:
\[ (\PMS_2^*) \ \ \   \ \ \ \ \ \ \ \ \ \ \ \ \ \ \ \ \   \ \ \ \ \ \ \ \ \ \ \ \ \ \ \ \
\forall_{i,j \in I}\big( \mu(i) = \mu(i\wedge j) + \mu(i \sim j)\big). \ \ \ \ \ \ \ \ \ \ \ \ \ \ \ \ \ \ \ \ \ \ \ \ \ \ \ \ \ \ \ \ \ \ \ \ \ \ \]
\end{defi}

\begin{rem}
Condition $(\PMS_2^*)$ involves quantification over a set and is absolutely
safe from a predicative point of view. Actually, it is only $(\PMS_2^*)$ that is needed
to construct the pre-integration space of simple functions.
\end{rem}

\begin{cor}\label{cor: premeasure1}
Let $\C M (\B \Lambda) := (X, I, \mu)$ be a pre-measure space and $i, j \in I$.
\begin{enumerate}[label=(\roman*)]
\item The operations $\vee$, $\wedge$ and $\sim$ are functions, and the triplet $(I, \vee, \wedge)$
satisfies the properties of a distributive lattice, except from
the absorption equalities.
\item $i \leq j \TOT
\B \lambda_0 (i) \subseteq \B \lambda_0 (j)$.
\end{enumerate}
\end{cor}

\begin{proof}
(i) We show that $\vee$ is a function, and for $\wedge$ and $\sim$ we proceed similarly. We have that
\begin{align*}
i = i{'} \ \& \ j = j{'} & \To \B \lambda_0 (i) = \B \lambda_0 (i{'}) \ \& \ \B \lambda_0 (j) = \B \lambda_0 (j{'})\\
& \To \B \lambda_0 (i) \vee \B \lambda_0 (j) = \B \lambda_0 (i{'}) \vee \B \lambda_0 (j{'})\\
& \To \B \lambda_0 (i \vee j) = \B \lambda_0 (i{'} \vee j{'})\\
& \To i \vee j = i{'} \vee j{'}.
\end{align*}
All properties
follow from the corresponding properties of complemented subsets (Proposition~\ref{prp: second2}),
from $(\PMS_1)$, and the fact that $\B \Lambda(X) \in \Set(I, \B X)$.
E.g., to show $i \vee j = j \vee i$,
we use the equalities $\B \lambda_0 (i \vee j) = \B \lambda_0 (i) \vee \B \lambda_0 (j) = \B \lambda_0 (j)
\vee \B \lambda_0 (i) = \B \lambda_0 (j \vee i)$. For the rest of the proof we proceed similarly.
\end{proof}

\noindent 
Next we give a constructive treatment of the classical Dirac measure as a pre-measure
on a set of integrable indices $I$. First we consider the total case, where $I := \D F(X, \D 2)$
is a Boolean algebra.

\begin{prop}\label{prp: Dirac1}

Let $\big(X, =_X, \neq_{(X, \mathsmaller{\D F(X, \D 2)})}\big)$ be a set inhabited by some $x_0 \in X$, and let
the maps $\vee, \wedge, \sim : \D F(X, \D2) \times \D F(X, \D2) \to \D F(X, \D2)$ and $\sim \colon \D F(X, \D2)
\to \D F(X, \D2)$, defined by the corresponding rules given in Definition~\ref{def: oppartial} for partial functions.
If $\B \Delta(X) := \big(\delta_0^1, \C E^{1}, \delta_1^1, \delta_0^0, \C E^{0}, \delta_1^0\big) $
is the set of complemented detachable subsets of $X$, where
by Example~\ref{ex: detachablefam}
\[ \B \delta_0(f) := \big(\delta_0^1(f), \delta_0^0(f)\big) := \big([f = 1], [f = 0]\big),\]
and if $\mu_{x_0} \colon \D F (X, \D 2) \sto [0, + \infty)$ is defined by the rule
\[ \mu_{x_0} (f) := f(x_0) =: \chi_{\B \delta_0(f)}(x_0); \ \ \ \ f \in \D F(X, \D 2), \]
then the triplet $\C M(\B \Delta(X)) := (X, \D F(X, \D 2), \mu_{x_0})$ is a pre-measure space.
\end{prop}

 \begin{proof}
 Straightforward calculations as in the proof of Proposition~\ref{prp: second1}(iii)
 prove the required equalities between complemented in
 condition $(\PMS_1)$. Clearly, the operation $\mu_{x_0}$ is a
 function, and a simple case-distinction shows the required equality
 $f(x_0) + g(x_0) = [f(x_0) \vee g(x_0)] + [f(x_0) \wedge g(x_0)]$.
 Let $f \in \D F (X, \D 2)$ and $\B B := (B^1, B^0) \in \C P^{\Disj}(X)$
 with $\B \alpha_0 (0) := \B B$. If $g \in \D F(X, \D 2)$ such that
 \begin{align*}
   \B \delta_0(f) \wedge \B B &:= \big([f = 1] \cap B^1, ([f=1] \cap B^0) \cup ([f=0] \cap B^1) \cup ([f=0] \cap B^0)\big) \\ &=
 \big([g=1], [g=0]\big),
 \end{align*}
 then the equality between the following complemented subsets
 \[\B \delta_0(f) - \B B = \big([f = 1] \cap B^1, ([f=1] \cap B^1) \cup ([f=0] \cap B^0) \cup ([f=0]
 \cap B^1)\big),\]
 \[\B \delta_0(f \sim g) = \big([f=1] \cap [g=0], [f=1=g] \cup [f=0=g] \cup [f=0] \cap [g=1]\big)\]
follows after considering all necessary cases (the rule Ex falsum quodlibet is necessary to this proof).
The required equality $f(x_0) = g(x_0) + f(x_0) \wedge (1 - g)(x_0))$ follows after considering all cases.
As $\mu_{x_0}(1) = 1 > 0$, $(\PMS_3)$ follows. For the proof of $(\PMS_4)$ we fix $\alpha : \Nat \to \D F(X, \D 2)$,
and we suppose that
 \begin{align*}
& \exists  \lim_{\mathsmaller{m \to + \infty}} \mu_{x_0}\bigg(\bigwedge_{n = 0}^m \alpha(n)\bigg) \  \& \
\lim_{\mathsmaller{m \to + \infty}} \mu_{x_0}\bigg(\bigwedge_{n = 0}^m \alpha(n)\bigg) > 0 \ \TOT \\
& \exists  \lim_{\mathsmaller{m \to + \infty}} \bigg(\bigwedge_{n = 0}^m \alpha(n)\bigg)(x_0) \ \& \
\lim_{\mathsmaller{m \to + \infty}} \bigg(\bigwedge_{n = 0}^m \alpha(n)\bigg)(x_0) > 0 \ \TOT \\
 & \exists  \lim_{\mathsmaller{m \to + \infty}} \prod_{n = 0}^m [\alpha(n)] (x_0)  \ \& \
\lim_{\mathsmaller{m \to + \infty}} \prod_{n = 0}^m [\alpha(n)](x_0) > 0.
 \end{align*}
Finally, we have that
\begin{align*}
 \lim_{\mathsmaller{m \to + \infty}} \prod_{n = 0}^m [\alpha(n)](x_0) > 0 & \To
 \lim_{\mathsmaller{m \to + \infty}} \prod_{n = 0}^m [\alpha(n)](x_0) =1\\
 & \TOT \exists_{m_0 \in \Nat}\forall_{m \geq m_0}\bigg(\prod_{n = 0}^m [\alpha(n)](x_0) = 1\bigg)\\
& \To \forall_{n \in \Nat}\big([\alpha(n)] (x_0) = 1\big)\\
 & \TOT x_0  \in \bigcap_{n \in \Nat}\delta_0^1 (\alpha(n)).\qedhere
 \end{align*}
\end{proof}

\begin{rem}
    Although the derivation of $(\PMS_2)$ in the above proof requires the Ex falsum
    quodlibet rule, the derivation of $(\PMS_2^*)$ rests on trivial calculations.
    Hence the whole proof in the latter case can be carried out in minimal logic!
\end{rem}

If partial functions are considered, then using Proposition~\ref{prp: second1}(iii)-(iv) we get similarly the following constructive version of the Dirac measure.

\begin{prop}\label{prp: Dirac2}
Let $\Lambda(X, \D 2) := (\lambda_0, \C E, \lambda_1, \F f)$ be a family of $($strongly extensional$)$
partial functions from $X$ to $\D 2$, with $\B f_i := (\lambda_0(i), \C E_i, \F f_i) \in \F F^{\se}(X, \D 2)$,
for every $i \in I$. Let $\vee, \wedge, \sim$ be operations on $I$, such that for every $i, j \in I$ we have that
$\B f_{i \vee j} = \B f_i \vee \B f_j,  \B f_{i \wedge j} = \B f_i \wedge \B f_j$, and
$\B f_{i \sim j} = \B f_i \sim \B f_j$. Moreover, let $\sim$ be an operation on $I$, such that
if $\B B$ is a given complemented subset of $X$, then the equality
$\delta(\B f_i) \wedge \B B = \delta(\B f_k)$ implies $\delta(\B f_i) - \B B = \delta(\B f_{i \sim k})$,
where the assignment routine $\delta$ is defined in Proposition~\ref{prp: compl3}.
If $x_0 \in X$ such that $x_0 \in \bigcap_{i \in I}\lambda_0(i),$
and if $\mu_{x_0} \colon I \sto [0, + \infty)$ is defined by the rule $ \mu_{x_0} (i) := \F f_i(x_0)$,
for every $i \in I$,
then the triplet $\C M(\B \Delta(I, X)) := (X, I, \mu_{x_0})$ is a pre-measure space.

\end{prop}

\section{Real-valued, partial functions}
\label{sec: realpartial}

Next we present the facts on real-valued, partial functions that are necessary for the definition of an
integration space within $\BST$ (Definition~\ref{def: preintspace}).

\begin{defi}\label{def: realpfunction}
If $(X, =_X, \neq_X)$ is an inhabited set, let $\B f_A := (A, i_A, f_A) \in \F F(X, \Real)$.
We call $\B f_A$ strongly extensional, if
$f_A$ is strongly extensional, where $A$ is equipped with its canonical inequality as a
subset of $X$ i.e., $f_A(a) \neq_{\Real} f_A(a{'}) \To i_A(a) \neq_X i_A(a{'})$,
for every $a, a{'} \in A$ $($where $a \neq_{\Real} b :\TOT a < b \vee b < a$, for every $a, b \in \Real)$.
Let $\F F(X) := \F F(X, \Real)$ be the class of partial functions
from $X$ to $\Real$, and $\F F^{\se}(X)$ the class of strongly extensional,
partial functions from $(X =_X, \neq_X)$ to $(\Real, =_{\Real}, \neq_{\Real})$.
Let $|\B f_A| := (A, i_A, |f_A|)$. If $\B f_B := (B, i_B, f_B)$ in $\F F(X)$ and $\lambda \in \Real$
\begin{center}
\begin{tikzpicture}

\node (E) at (0,0) {$A$};
\node[right=of E] (B) {$X$};
\node[right=of B] (F) {$B$};
\node[below=of B] (C) {$\Real$,};

\draw[right hook->] (E)--(B) node [midway,above] {$i_A$};
\draw[left hook->] (F)--(B) node [midway,above] {$i_B$};
\draw[->] (E)--(C) node [midway,left] {$f_A \ \ $};
\draw[->] (F)--(C) node [midway,right] {$ \ f_B$};

\end{tikzpicture}
\end{center}
let $\lambda\B f_A := (A, i_A, \lambda f_A) \in \F F(X)$ and
$\B f_A \ \square \ \B f_B := \big(A \cap B, i_{A \cap B}, \big(f_A \ \square
\ f_B\big)_{A \cap B}\big)$, where
\[ (f_A \ \square \ f_B\big)_{A \cap B}(a, b) := f_A(a) \
\square \ f_B(b); \ \ \ \ (a, b) \in A \cap B, \  \ \square \in \{+, \cdot, \wedge, \vee\}.\]
The totality of $I$-families of strongly extensional, partial functions i.e., of structures
$\Lambda(X,\Real) := (\lambda_0, \C E, \lambda_1, \F f)$, with $\B f_i := (\lambda_0(i), \C E_i, \F f_i)$
strongly extensional, for every $i \in I$, is denoted by $\Fam^{\se}(I, X, \Real)$.
\end{defi}

The operation $(f_A \ \square \ f_B\big)_{A \cap B} \colon A \cap B \sto \Real$
is a function. If $(a, b) =_{A \cap B} (a{'}, b{'}) :\TOT i_A =_X i_A(a{'}) \TOT a =_A a{'}$,
we get $f_A(a) =_{\Real} f_A(a{'})$.
Since $i_B(b) =_X i_A(a)$ and $i_B(b{'}) =_X i_A(a{'})$,
we also get $b =_B b{'}$ and hence $f_B(b) =_{\Real} f_B(b{'})$.
If $\lambda$ denotes also the constant function $\lambda \in \Real$ on $X$
we get as a special case the partial function
$\B f_A \wedge \B \lambda := \big(A \cap X, i_{A \cap X}, \big(f_A \wedge \lambda\big)_{A \cap X}\big)$,
where $A \cap X := \{ (a, x) \in A \times X\mid i_A(a) =_X x\}$, $i_{A \cap X}(a, x) := i_A(a)$,
and $\big(f_A \wedge \lambda\big)_{A \cap X}(a,x):= f_A(a) \wedge \lambda(x) := f_A(a) \wedge \lambda$,
for every $(a, x) \in A \cap X$.
By Definition~\ref{def: famofpartial}, if $\Lambda(X,\Real) := (\lambda_0, \C E, \lambda_1, \F f)
\in \Fam(I, X, \Real)$, then $\B f_i := \big(\lambda_0(i), \C E_i, \F f_i\big) \in \F F(X)$, for every $i \in I$,
and if $i =_I j$, the following diagrams commute
\begin{center}
\begin{tikzpicture}

\node (E) at (0,0) {$\lambda_0(i)$};
\node[right=of E] (B) {};
\node[right=of B] (F) {$\lambda_0(j)$};
\node[below=of B] (A) {$X$};
\node[below=of A] (C) {$ \ \Real$.};

\draw[->,bend left] (E) to node [midway,below] {$\lambda_{ij}$} (F);
\draw[->,bend right=50] (E) to node [midway,left] {$\F f_i$} (C);
\draw[->,bend left=50] (F) to node [midway,right] {$\F f_j$} (C);
\draw[right hook->,bend right=20] (E) to node [midway,left] {$ \ \C E_i \ $} (A);
\draw[left hook->,bend left=20] (F) to node [midway,right] {$\ \C E_j \ $} (A);

\end{tikzpicture}
\end{center}
If $\B f_i$ is strongly extensional, then, for every $u, w \in \lambda_0(i)$, we get
$\F f_i(u) \neq_{\Real} \F f_i(w) \To \C E_i(u) \neq_X \C E_i(w)$.
We may also regard $\Lambda(X,\Real)$ as a $\lambda_0I (X, \Real)$-set
of real valued, strongly extensional, partial functions,
following the construction in Definition~\ref{def: setofpartialfunctions}.

\begin{defi}\label{def: subseq}
Let $\Lambda(X,\Real) := (\lambda_0, \C E, \lambda_1, \F f) \in \Fam^{\se}(I, X, \Real)$.
We write $g \colon \lambda_0I (X, \Real) \to Y$ to denote a function $g \colon I \to Y$, where
$I$ is equipped with the equality in Definition~\ref{def: setofpartialfunctions},
and we may also write $g(\B f_i)$ instead of $g(i)$.
If $\kappa \colon \Nat^+ \to I$, the family
$$\Lambda(X, \Real) \circ \B \kappa := \big(\lambda_0 \circ \kappa, \C E \circ \kappa, \lambda_1
\circ \kappa, \F f \circ \kappa\big) \in \Fam(\Nat^+, X, \Real)$$
is the $\kappa$-subsequence
of $\Lambda(X, \Real)$, where
$(\lambda_0 \circ \kappa)(n) := \lambda_0(\kappa(n))$, $\big(\C E \circ \kappa\big)_n :=\C E_{\kappa(n)}$,
$(\lambda_1 \circ \kappa)(n, n) := \lambda_{\kappa(n)\kappa(n)} :=\id_{\lambda_0(\kappa(n))}$ and
$\big(\F f \circ \kappa\big)_n := \F f_{\kappa(n)}$
for every $n \in \Nat^+$.
\end{defi}

If we consider the intersection
$\bigcap_{n \in \Nat^+}(\lambda_0 \circ \kappa)(n) :=\bigcap_{n \in \Nat^+}\lambda_0(\kappa(n))$,
by Definition~\ref{def: intfamilyofsubsets} we get
\[ \Phi \in \bigcap_{n \in \Nat^+}\lambda_0(\kappa(n))  :\TOT \Phi \colon
\bigcurlywedge_{n \in \Nat^+}\lambda_0(\kappa(n)) \ \ \& \ \
\forall_{n, m \in \Nat^+}\big(\C E_{\kappa(n)}(\Phi_n) =_X \C E_{\kappa(m)}(\Phi_m)\big), \]
\[ \Phi =_{\mathsmaller{\bigcap_{n \in \Nat^+}\lambda_0(\kappa(n))}} \Theta
:\TOT \C E_{\kappa(1)}(\Phi_1) =_X \C E_{\kappa(1)}(\Theta_1), \]
\[ e^{\Lambda(X, \Real) \circ \kappa} \colon
\bigcap_{n \in \Nat^+}\lambda_0(\kappa(n)) \eto X, \ \ \ \
e^{\Lambda(X, \Real) \circ \kappa}(\Phi) :=
\big(\C E \circ \kappa)_1 (\Phi_1) := \C E_{\kappa(1)}(\Phi_1).\]

\begin{defi}\label{def: sigmapfun}
Let $\Lambda(X,\Real) := (\lambda_0, \C E, \lambda_1, \F f) \in
\Fam(I, X, \Real)$, $\kappa \colon \Nat^+ \to I$, and
$\Lambda(X, \Real) \circ \kappa$ the $\kappa$-subsequence of $\Lambda(X, \Real)$.
If $(A, i_A) \subseteq \bigcap_{n \in \Nat^+}\lambda_0(\kappa(n))$, we define the function
\[ \sum_{n \in \Nat^+} \F f_{\kappa(n)} \colon A \to \Real, \ \ \ \
\bigg(\sum_{n \in \Nat^+}\F f_{\kappa(n)}\bigg)(a)
:= \sum_{n \in \Nat^+}\F f_{\kappa(n)}\bigg(\big[i_A(a)\big]_n\bigg); \ \ \ \ a \in A,\]
under the assumption that the series on the right converges in $\Real$, for every $a \in A$.
\end{defi}

In the special case
$\big(\bigcap_{n \in \Nat^+}\lambda_0(\kappa(n)), \id_{\mathsmaller{\bigcap_{n \in \Nat^+}\lambda_0(\kappa(n))}}\big)
\subseteq \bigcap_{n \in \Nat^+}\lambda_0(\kappa(n))$, we get the function
\[ \sum_{n \in \Nat^+}\F f_{\kappa(n)} \colon
\bigcap_{n \in \Nat^+}\lambda_0(\kappa(n)) \to \Real, \ \ \ \
\bigg(\sum_{n \in \Nat^+}\F f_{\kappa(n)}\bigg)(\Phi)
:= \sum_{n \in \Nat^+}\F f_{\kappa(n)}\big(\Phi_n\big); \ \ \ \
\Phi \in \bigcap_{n \in \Nat^+}\lambda_0(\kappa(n)),\]
under the same convergence assumption. The following fact is shown in~\cite{Pe20},
pp.~212--213, and it is used in Definition~\ref{def: canonical}.

\begin{prop}\label{prp: sextpf1}
If in Definition~\ref{def: sigmapfun} the partial functions $\B f_{\kappa(n)} :=
\big(\lambda_0(\kappa(n)), \C E_{\kappa(n)}, \F f_{\kappa(n)}\big)$ are strongly extensional,
for every $n \in \Nat^+$, then the real-valued, partial function
\[ \B f_A := \bigg(A, \ e^{\Lambda(X, \Real) \circ \kappa}
\circ i_A, \ \sum_{n \in \Nat^+} \F f_{\kappa(n)} \bigg)\]
\begin{center}
\begin{tikzpicture}

\node (E) at (0,0) {$A$};
\node[right=of E] (G) {$\bigcap_{n \in \Nat^+}\lambda_0(\kappa(n))$};
\node[right=of G] (H) {};
\node[right=of H] (F) {$X$};
\node[below=of F] (A) {$\Real$};

\draw[right hook->] (E)--(G) node [midway,above] {$i_A$};
\draw[right hook->] (G)--(F) node [midway,above] {$e^{\Lambda(X, \Real) \circ \kappa} $};
\draw[right hook->,bend right=20] (E) to node [midway,left] {$
\sum_{n \in \Nat^+} \F f_{\kappa(n)} \ \ \ \  $} (A);

\end{tikzpicture}
\end{center}
is strongly extensional.
\end{prop}

\section{Pre-integration spaces}
\label{sec: preint}

In this section, and in accordance to our previous predicative
reconstruction of Bishop-Cheng measure space, we introduce the notion
of a pre-integration space as a predicative counterpart to the notion
of an integration space in $\BCMT$. The notion of a Bishop-Cheng
integration space is defined in~\cite{BB85}, p.~217, and appeared
first in~\cite{BC72}, p.~2. Condition $(\BCIS_2)$ is the constructive
counterpart to Daniell's classical continuity condition in the
definition of a Daniell space.  The exact relation of a Bishop-Cheng
integration space to that of a Daniell space is explained in~\cite{Pe22c}.

\begin{defiC}[\textbf{Bishop-Cheng integration space}]
A triplet $(X, L, \int)$ is a (Bishop-Cheng) integration space if $(X, =_X, \neq_X)$ is an inhabited set with inequality,
$L$ is a subset of $\F F^{\se}(X)$, and $\int \colon L \to \Real$,
such that the following properties hold\index{Bishop-Cheng definition of an integration space}.
\begin{enumerate}[label=$(\BCIS_{\arabic*})$]
\item If $\B f, \B g \in L$ and $\alpha, \beta \in \Real$, then $\alpha \B f + \beta \B g $, $|\B f|$, and
$\B f \wedge \B 1$ belong to $L$, and
\begin{align*}
  \int(\alpha \B f + \beta \B g) = \alpha \int \B f + \beta \int \B g.
\end{align*}
\item If $\B f \in L$ and $(\B f_n)$ is a sequence of non-negative functions in $L$ such that
$\sum_{n}\int \B f_n$ converges and $\sum_n \int(\B f_n) < \int \B f$, then there exists $x \in X$
such that $\sum_n \F f_n(x)$ converges and $\sum_n \F f_n(x) < \F f(x)$.
\item There exists a function $\B p$ in $L$ with $\int \B p = 1$.
\item For each $\B f \in L$, $\lim_{n \to \infty} \int(\B f \wedge \B n) = \int \B f$ and
$\lim_{n \to \infty} \int(|\B f| \wedge \B n^{-1}) = 0$.
\end{enumerate}
\end{defiC}
\noindent
As already mentioned in the introduction, there is no explanation how the set $L$ is ``separated''
from the proper class
$\F F^{\se}(X)$, so that the integral $\int$ can be defined as a real-valued function on $L$.
The extensional character of $L$ is also not addressed.
This impredicative approach to $L$ is behind the simplicity of the Bishop-Cheng integration space.
E.g., in condition $(\BCIS_4)$ the formulation of the limit is immediate as the terms $\B f \wedge \B n \in L$
and $\int$ is defined on $L$.
If one predicatively reformulates the Bishop-Cheng definition though, where
an $I$-family of strongly extensional, real-valued, partial functions
is going to be used instead of $L$, then one needs to use an element
$\alpha(n)$ of the index-set $I$ such that $\B f \wedge \B n =_{\F F^{\se}(X)} \B f_{\alpha(n)}$,
in order to express the corresponding
limit. The formulation of the continuity condition $(\BCIS_2)$ takes
the form
 \[ \  \ \ \ \  \ \ \ \ \ \ \ \ \  \ \ \ \ \ \ \ \  \forall_{i \in I}\forall_{\kappa
 \in \D F(\Nat^+, I)}\bigg\{\bigg[\sum_{n \in \Nat^+}\int \B f_{\kappa(n)} \in \Real \  \ \& \
 \sum_{n \in \Nat^+}\int \B f_{\kappa(n)} < \int \B f_i\bigg] \To  \ \ \ \ \ \ \ \ \ \ \ \ \ \ \ \ \ \ \ \ \
 \ \ \ \ \ \ \ \]
 \[ \ \  \exists_{(\Phi, u) \ \in \ \big(\bigcap_{n \in \Nat^+}\lambda_0(\kappa(n)) \big)\ \cap \
 \lambda_0(i)}\bigg( \sum_{n \in \Nat^+}\F f_{\kappa(n)}(\Phi_n) \in \Real \  \ \& \ \
 \sum_{n \in \Nat^+}\F f_{\kappa(n)}(\Phi_n) < \F f_i(u) \bigg)\bigg\},\]
where
\[ \bigg(\bigcap_{n \in \Nat^+}\lambda_0(\kappa(n))\bigg) \cap  \lambda_0(i)
:= \bigg\{(\Phi, u) \in \bigg(\bigcap_{n \in \Nat^+}\lambda_0(\kappa(n))\bigg) \times
\lambda_0(i) \ \mathlarger{\mathlarger{\mid}} \ \C E_{\kappa(1)}(\Phi) =_X \C E_i(u)\bigg\},\]
\begin{center}
\begin{tikzpicture}

\node (E) at (0,0) {$\lambda_0(\kappa(n))$};
\node[right=of E] (B) {$X$};
\node[right=of B] (F) {$\lambda_0(i)$};
\node[below=of B] (C) {$\Real$.};

\draw[right hook->] (E)--(B) node [midway,above] {$\C E_{\kappa(n)}$};
\draw[left hook->] (F)--(B) node [midway,above] {$\C E_i$};
\draw[->] (E)--(C) node [midway,left] {$\F f_{\kappa(n)} \ \ $};
\draw[->] (F)--(C) node [midway,right] {$ \ \F f_i$};

\end{tikzpicture}
\end{center}
Next we directly formulate the Bishop-Cheng definition of an
integration space using appropriate operations on the index-set of the
appropriate family (set) of real-valued partial functions that
replaces the original impredicative subset (actually, proper-class)
$L$ of $\F F^{\se}(X)$.

\begin{defi}[Pre-integration space within $\BST$]\label{def: preintspace}
Let $(X, =_X, \neq_X)$ be an inhabited set with inequality, and let the set $(I, =_I)$ be equipped with
operations $\cdot_a \colon I \sto I$, for every $a \in \Real$, $+ \colon I \times I \sto I$, $|.| \colon I \sto I$,
and $\wedge_1 \colon I \sto I$, where
\[ \cdot_a(i) := a \cdot i, \ \ \ \ +(i,j) := i + j, \ \ \ \ |.|(i) := |i|; \ \ \ \ i \in I, \ a \in \Real.\]
Let also the operation $\wedge_a \colon I \sto I$, defined by the previous operations through the rule
\[ \wedge_a := \cdot_a \circ \wedge_1 \circ \cdot_{a^{-1}}; \ \ \ \ a \in \Real \ \& \ a > 0.\]
Let $\Lambda(X,\Real) := (\lambda_0, \C E, \lambda_1, \F f) \in \Set^\se(I, X, \Real)$ i.e.,
$\B f_i =_{\F F^{\se}(X)} \B f_j \To i =_I j$, for every $i, j \in I$, and $\B f_i := \big(\lambda_0(i),
\C E_i, \F f_i\big)$ is strongly extensional, for every $i \in I$. Let $\int \colon I \to \Real$ be a function, where $i \mapsto \int i$, for every $i \in I$,
such that the following conditions hold:
\[(\PIS_1) \ \  \ \ \ \ \ \ \ \ \ \ \ \ \ \ \ \ \ \ \ \  \ \ \ \ \ \ \
\forall_{i \in I}\forall_{a \in \Real}\bigg(a \B f_i = \B  f_{a \cdot i} \ \ \& \ \
\int a \cdot i = a \int i\bigg). \ \ \ \ \ \ \ \ \ \ \ \ \ \ \ \ \ \ \ \ \ \
\ \ \ \ \ \ \ \ \ \ \ \ \ \ \ \ \ \ \  \]
\[(\PIS_2) \ \  \ \ \ \ \ \ \ \ \ \ \ \ \ \ \ \ \ \ \ \  \ \
\forall_{i, j \in I}\bigg(\B f_i + \B f_j = \B f_{i + j} \  \ \& \ \
\int (i + j) = \int i \ + \int j\bigg). \ \ \ \ \ \ \ \ \ \ \ \ \ \ \ \ \ \ \ \ \ \ \ \ \ \ \ \ \ \ \ \ \ \ \  \]
\[(\PIS_3) \ \  \ \ \ \ \ \ \ \ \ \ \ \ \ \ \ \ \ \ \ \ \ \ \ \ \ \ \ \ \ \ \ \ \ \ \ \
\ \  \ \ \forall_{i \in I}\big( |\B f_i| = \B f_{|i|} \big). \ \ \ \ \ \ \ \ \ \ \ \ \ \ \ \ \ \ \ \ \ \
\ \ \ \ \ \ \ \ \ \ \ \ \ \ \ \ \ \ \ \ \ \ \ \ \ \ \ \ \ \ \ \ \ \   \]
\[(\PIS_4) \ \  \ \ \ \ \ \ \ \ \ \ \ \ \ \ \ \ \ \ \ \ \ \ \ \ \ \ \ \ \ \ \ \ \ \ \ \
\ \ \  \ \ \forall_{i \in I}\big( \B f_i \wedge \B 1 = \B f_{\wedge_1(i)} \big). \ \
\ \ \ \ \ \ \ \ \ \ \ \ \ \ \ \ \ \ \ \ \ \ \ \ \ \ \ \ \ \ \ \ \ \ \ \ \ \ \ \ \ \ \ \ \ \ \ \ \ \ \ \ \   \]
\[(\PIS_5) \  \ \ \ \  \ \ \ \ \ \ \ \ \  \ \ \ \ \ \ \ \  \ \
\forall_{i \in I}\forall_{\kappa \in \D F(\Nat^+, I)}\bigg\{\bigg[\sum_{n \in \Nat^+}\int
\kappa(n) \in \Real \  \ \& \ \sum_{n \in \Nat^+}\int \kappa(n) < \int i\bigg]  \To  \ \ \
\ \ \ \ \ \ \ \ \ \ \ \ \ \ \ \ \ \ \ \ \ \ \ \ \  \]
\[ \exists_{(\Phi, u) \ \in \ \big(\bigcap_{n \in \Nat^+}\lambda_0(\kappa(n)) \big)\ \cap \
\lambda_0(i)}\bigg( \sum_{n \in \Nat^+}\F f_{\kappa(n)}(\Phi_n) \in \Real  \ \& \
\sum_{n \in \Nat^+}\F f_{\kappa(n)}(\Phi_n)
< \F f_i(u). \bigg)\bigg\}\]
\[(\PIS_6) \ \  \ \ \ \ \ \ \ \ \ \ \ \ \ \ \ \ \ \ \ \ \ \ \ \ \ \ \ \ \ \ \ \ \ \ \ \  \
\ \ \ \ \ \  \ \ \exists_{i \in I}\bigg( \int  i =_{\Real} 1 \bigg). \ \ \ \ \ \ \  \ \ \
\ \ \ \ \ \ \ \ \ \ \ \ \ \ \ \ \ \ \ \ \ \ \ \ \ \ \ \ \ \ \ \ \ \ \ \ \ \ \ \ \ \ \ \   \]
\[(\PIS_7) \ \  \ \ \ \ \ \ \ \ \ \ \ \ \ \ \ \ \ \ \ \ \ \ \ \ \
\ \   \forall_{i \in I}\bigg(\lim_{\mathsmaller{n \longrightarrow +\infty}}
\int \wedge_n(i) \in \Real \ \ \& \ \lim_{\mathsmaller{n \longrightarrow +\infty}}
\int \wedge_n(i) = \int i\bigg). \  \ \ \ \ \ \ \ \ \ \ \  \ \ \ \ \ \ \ \ \ \ \ \ \ \ \ \ \ \ \]
\[(\PIS_8) \ \  \ \ \ \ \ \ \ \ \ \ \ \ \ \ \ \ \ \ \ \ \ \ \ \ \
\ \   \forall_{i \in I}\bigg(\lim_{\mathsmaller{n \longrightarrow +\infty}}
\int \wedge_{\frac{1}{n}}(|i|) \in \Real \ \ \& \
\lim_{\mathsmaller{n \longrightarrow +\infty}} \int \wedge_{\frac{1}{n}}(|i|) = 0\bigg). \ \
\ \ \ \ \ \ \ \ \ \ \ \ \ \ \ \ \ \ \ \ \ \ \ \ \   \]
We call the structure $\C L_0 := \big(X, I, \Lambda(X,\Real), \int\big)$ a pre-integration
space\index{pre-integration space}.

\end{defi}

All the operations on $I$ defined above are functions. E.g., since $\Lambda(X, \Real)  \in \Set^\se(I, X, \Real)$,
\[ i =_I i{'} \To  \B f_i =  \B f_{i{'}} \To a  \B f_i = a  \B f_{i{'}} \To  \B f_{a \cdot i}
=  \B f_{a \cdot i{'}} \To a \cdot i =_I a \cdot i{'}.\]
The most fundamental example of an integration space within $\BCMT$ is
that induced by a positive measure $\mu$ on a locally compact metric
space $X$ i.e., a non-zero linear map on the functions with compact
support $C^{\supp}(X)$ (see~\cite{BB85}, pp.~220-221). In~\cite{Gr22,GP22}
this major example is described as a pre-integration
space. For that a notion of a locally compact metric space with a
modulus of local compactness is introduced. If $(X, d)$ is an
inhabited metric space with $x_0 \in X$, and $(K_n)_{n \in \Nat}$ is a
sequence of compact subsets of $X$, a \textit{modulus of local compactness} for $X$ is a function
$\B \kappa \colon \Nat \to \Nat$,
$n \mapsto \B \kappa(n),$
such that $[d_{x_0} \leq n] \subseteq K_{\B \kappa(n)}$, for every $n
\in \Nat$, where $d_{x_0} \colon X \to [0, + \infty)$ is defined by
$d_{x_0}(x) := d(x, x_0)$, for every $x \in X$.  In this way the
initial impredicativity of Bishop's notion of a locally compact
metric space (for \textit{every} bounded subset $B$ of $X$,
\textit{there is} a compact subset $K$ of $X$ with $B \subseteq K$)
is avoided.  If $\big(X, d, (K_n)_{n \in \Nat}, \B \kappa\big)$ is a
locally compact metric space with a modulus of local compactness, a
uniformly continuous function on every bounded subset of $X$ (this
impredicativity can be easily avoided) has compact support if there
is $m \in \Nat$ such that $K_m$ is a support of $f$ i.e.,
$\forall_{x \in X}\big(d(x, K_m) > 0 \To f(x) = 0\big)$.  If we
consider their set $C^{\supp}(X)$ as the index-set of the family
$\Supp(X, \Real) \in \Set^\se(I, X, \Real)$ of strongly extensional,
real-valued, partial functions on $X$
$$f \mapsto (X, \id_X, f),$$
where $X$ is equipped the canonical
inequality induced by its metric $\big(x \neq_{(X,d)} x{'} :\TOT d(x, x{'})> 0\big)$,
then the following result is shown
in~\cite{Gr22, GP22} within $\BST$, and it is the starting point of
a predicative reconstruction of the integration theory of locally
compact metric spaces within $\BST$.

\begin{thm}[The pre-integration space of a locally compact metric space with a modulus of local
compactness]\label{thm: intlcms}
Let $\big(X, d, (K_n)_{n \in \Nat}, \B \kappa\big)$ be a locally compact metric space, $\neq_{(X,d)}$ the
canonical inequality on $X$, and let $I := C^{\supp}(X)$ be equipped with the following operations:
\begin{enumerate}[label=(\roman*)]
\item If $a \in \Real$, then
$\cdot_a \colon I \to I$ is defined by $f \mapsto af$.
\item $+ \colon I \times I \to I$ is the addition of functions on $I$.
\item $|.| \colon I \to I$ is defined by $f \mapsto |f|$.
\item $\wedge_1 \colon I \to I$ is defined by $f \mapsto f \wedge 1$.
\item If $a > 0 \in \Real$, then $\wedge_a \colon I \to I$ is defined as the composition
$\wedge_a := \cdot_a \circ \wedge_1 \circ \cdot_{a^{-1}}$.
\end{enumerate}
Let the obviously defined set $\Supp(X, \Real)$ of strongly extensional
real-valued, partial functions over $I$. If $\mu \colon I \to \Real$ is a linear, positive measure on $X$ i.e.,
there is $f \in I$ with $\mu(f) > 0$, and for every $f \in I$ we have that
$f \geq 0 \To \mu(f) \geq 0$, let
\[ \int\_ d\mu \colon I \to \Real, \ \ \ \ f \mapsto \int f d\mu := \mu(f); \ \ \ \ f \in I.\]
Then $\big(X, I, \Supp(X, \Real), \int\_ d\mu\big)$ is a pre-integration space.
\end{thm}

\section{Simple functions}
\label{sec: simple}

In this section we construct the pre-integration space of simple
functions from a given pre-measure space (Theorem~\ref{thm:
  pre-10.10}). This is a predicative translation within $\BST$ of the
construction of an integration space from the simple functions of a
measure space (Theorem 10.10 in~\cite{BB85}). Although we follow the
corresponding construction in section 10 of chapter 6 in~\cite{BB85}
closely, our approach allows us to not only work completely
predicatively, but also to carry out all proofs avoiding the axiom of
countable choice.
\textit{For the remainder of this section we fix an inhabited set with inequality
$(X,=_X,\neq_X)$, and an $I$-family
of complemented subsets} $\B \Lambda(X) := \big(\lambda_0^1, \C E^{1}, \lambda_1^1,
\lambda_0^0, \C E^{0},\lambda_1^0)$ with $i_0 \in I$. For every $i \in I$ let
\begin{align*}
  \B \chi_i := \big(\dom_i := \lambda_0^1(i) \cup \lambda_0^0(i), \C E_i, \chi_i\big) \in \F F^{\se}(X),
\end{align*}
where $\chi_i$ is the characteristic function of the complemented subset
$\B \lambda_0(i) := \big(\lambda_0^1(i),\lambda_0^0(i)\big)$ of $X$.

\begin{defi}\label{def: simple}
If $n \in \Nat^+$, $i_1, \ldots, i_n \in I$, and $a_{1}, \ldots a_{n} \in \Real$,
the triplet\footnote{The
  fact that $\sum_{k = 1}^n a_{k}\chi_{i_k}$ is strongly extensional is
  based on Remark~\ref{rem: compl2} and the following properties of $a,
  b \in \Real$: $a + b > 0 \To a > 0 \vee b > 0$ and $a \cdot b > 0 \To
  a \neq_{\Real} 0 \wedge b \neq_{\Real} 0$ (see~\cite{BB85}, p.~26
  and~\cite{Pe18}, p.~17, respectively).}
\[ \sum_{k = 1}^n a_{k} \B \chi_{i_k} := \bigg( \bigcap_{k = 1}^n \dom_{i_k},
i_{\mathsmaller{\bigcap_{k = 1}^n
\dom_{i_k}}}, \ \sum_{k = 1}^n a_{k}\chi_{i_k}\bigg) \in \F F^{\se}(X)\]
is called a simple function.
Consider the totality\footnote{The elements of $S(I,\B \Lambda(X))$ are
  pairs $(n,u)$, where $n\in\Nat^+$ and $u=\big((a_1, i_1), \ldots, (a_n, i_n)\big)$  is
  an $n$-tuple of pairs in $\Real\times I$. With a bit of abuse of notation we also write the
  elements of this Sigma-set as $(a_k,i_k)_{k=1}^n$,
  which is more convenient and contains all the information
  needed to write down the corresponding element in its proper form.}
$S(I, \B \Lambda(X)) := \sum_{n \in \Nat^+}(\Real \times I)^n$.
Let the non-dependent assignment routine $\dom_0 \colon  S(I, \B \Lambda(X)) \sto \D V_0$, defined by
$\dom_0 (n, u) := \bigcap_{k = 1}^n \dom_{i_k}$, for every $n \in \Nat^+$ and every
$u := \big((a_1, i_1), \ldots, (a_n, i_n)\big)$.
Furthermore, let $\C Z:\bigcurlywedge_{(n, u) \in S(I, \B \Lambda(X))}\mathbb{F}(\dom_0 (n, u),X)$ be the
dependent assignment routine, where $\C Z_{(n,u)}:\dom_0 (n, u)\hookrightarrow X$ is the canonical
embedding induced by the embeddings $\C E_{i_k}^{1}$ and $\C E_{i_k}^{0}$,
where $k \in \{1, \ldots, n\}$, and the dependent assignment routine
$\F f:\bigcurlywedge_{(n, u) \in S(I, \B \Lambda(X))}\mathbb{F}(\dom_0 (n, u),\Real)$ given by
$\F f_{(n,u)} := \sum_{k = 1}^n a_{k}\chi_{i_k}$, for every $n \in \Nat^+$ and every
$u := \big((a_1, i_1), \ldots, (a_n, i_n)\big)$.
We now take $S(I, \B \Lambda(X))$ to be equipped with the equality
\begin{align*}
  (a_k,i_k)_{k=1}^n=_{S(I, \B \Lambda(X))} (b_\ell,j_\ell)_{\ell=1}^m :\TOT \sum_{k = 1}^n a_k \B \chi_{i_k}=_{\F F^{\se}(X)} \sum_{\ell = 1}^m b_\ell \B \chi_{j_\ell}
\end{align*}
and define the set of \emph{simple functions} as the $S(I, \B \Lambda(X))$-set
of strongly extensional partial functions
$\Simple(\B \Lambda(X)) := \big(\dom_0, \C Z, \dom_1, \F f\big)\in \Set^{\se}\big(S(I, \B \Lambda(X)), X, \Real\big)$,
where $\dom_1$ is defined through dependent unique choice as explained in Definition~\ref{def: setofpartialfunctions}.
\end{defi}

It is immediate to show that $\Simple(\B \Lambda(X))$ is a set of
strongly extensional, partial functions over $S(I, \B
\Lambda(X))$. Next we translate the results from~\cite{BB85} needed to
prove that the simple functions form a pre-integration space. For the
most part the proofs of the many lemmas work exactly analogous to the
corresponding ones in~\cite{BB85}, so we won't give them here. Some of
the results can however be sharpened, thus allowing us to avoid the
axiom of countable choice altogether, and we present the proofs of
those results.  First, we state the predicative analogues of lemmas
(10.2) - (10.5) of chapter 6 in~\cite{BB85}.  One of the reasons
working always with an inhabited set $X$ is that an element $x_0 \in X$
is needed in Bishop's negativistic definition of the empty
subset\footnote{For a positively defined empty subset of $X$
see~\cite{PW22}.}  $\emptyset_X$ of $X$ (see~\cite{Bi67}, p.~65).

\begin{lem}\label{prp: pre-10.2}
Let $\C M (\B \Lambda) := (X, I, \mu)$ be a pre-measure space, $i_1,...,i_n\in I$, and
$F:=\bigcap_{k=1}^n\dom_{i_k}$.
\begin{enumerate}[label=(\roman*)]
\item If $i\in I$ such that $\lambda_0^1(i)= \emptyset_X$, then $\mu(i)=0$.
\item There is $j\in I$ such that $\bm\lambda_0(j)=(\emptyset_X,F)$.
\item If $j\in I$, then there is $k\in I$ such that
$\bm\lambda_0(k)=(\lambda_0^1(j)\cap F,\lambda_0^0(j)\cap F) $
  and $\mu(k)=\mu(j)$.
\item If $i,j\in I$, $F':= \dom_i \cap \dom_j \cap F$, and $\chi_i(x)\leq \chi_j(x)$, for every $x\in F'$,
  then $\mu(i)\leq \mu(j)$.
\end{enumerate}
\end{lem}

\begin{rem}
    For the proof of Lemma~\ref{prp: pre-10.2} it suffices to use
    condition $(\PMS_2^*)$. It is only here that we rely on $(\PMS_2)$,
    or $(\PMS_2^*)$.
\end{rem}

Many later proofs rely on the fact that we can restrict our attention to disjoint simple functions.
The next lemma, which corresponds to lemma (7.8) of chapter 6 in~\cite{BB85},
makes this fact precise.

\begin{lem}\label{prp: pre-7.8}\hfill
\begin{enumerate}[label=(\roman*)]
\item If $\overline{n} := \{1, \ldots, n\}$, the assignment routine $\disjrep
\colon S(I, \B \Lambda(X)) \sto S(I, \B \Lambda(X))$
\[ (a_k,i_k)_{k=1}^n \mapsto \bigg(\sum_{f(k)=1}a_k\; , \;j_f:=\bigg(\bigwedge_{f(k)=1}i_k\bigg)\sim\bigg(\bigvee_{f(k) =
0}i_k\bigg)\bigg)_{f \in \D F(\overline{n}, \D 2)}\]
is a function equal to $\id_{S(I, \B \Lambda(X))}$.
If $v\in S(I, \B \Lambda(X))$ and
$\disjrep(v):=(b_\ell,j_\ell)_{\ell=1}^m$, then the complemented subsets
$\bm\lambda_0(j_\ell)$ are \emph{disjoint}, i.e., if $\ell \neq k$, then
$\chi_{j_k}\cdot\chi_{j_\ell}=0$ on $\dom_{j_l} \cap \dom_{j_k}$.
\item For every $v:=(a_k,i_k)_{k=1}^n\in S(I, \B \Lambda(X))$ we have
\begin{align*}
\sum_{k=1}^n a_k\cdot\mu(i_k)=\sum_{f \in \D F(\overline{n}, \D 2)}\bigg(\sum_{f(k)=1}a_k\bigg)\cdot\mu(j_f)
\end{align*}
where $j_f\in I$ is defined, for every $f \in \D F(\overline{n}, \D 2)$, as in {\normalfont (i)}.
\end{enumerate}
\end{lem}

\begin{rem}\label{rem: simpleInduction}
If $v\in S(I, \B \Lambda(X))$, we call $\disjrep(v)$ the \emph{disjoint representation} of $v$.
The above lemma allows us to restrict our attention to disjoint simple functions,
whenever we want to prove a statement for all simple functions.
Combining this with the fact that we can prove statements about simple function by induction
on their first component $n\in\Nat^+$ we can show that some \textit{extensional} property $P$ holds for
all $v\in S(I, \B \Lambda(X))$ if can show the following:
\begin{itemize}
\item $P$ holds for all simple functions of length one $v :=(a,i)$, with $a\in\Real$ and $i\in I$.
\item If $P$ holds for $v=(a_k,i_k)_{k=1}^n$ disjoint
  and we have $a_{n+1}\in\Real$ and $i_{n+1}\in I$ disjoint from any of the indices
  $i_1,...,i_n\in I$, then $P$ also holds for $(a_k,i_k)_{k=1}^{n+1}\in S(I, \B \Lambda(X))$.
\end{itemize}
\end{rem}

\begin{lem}\label{lem: pre-10.6}
\hfill
\begin{enumerate}[label=(\roman*)]
\item Let $v:=(a_k,i_k)_{k=1}^n,w:=(b_\ell,j_\ell)_{\ell=1}^m\in S(I, \B \Lambda(X))$ such that
\begin{align*}
\sum_{k=1}^n a_k\cdot\chi_{i_k}(x)\leq\sum_{\ell=1}^m b_\ell\cdot\chi_{j_\ell}(x)
\end{align*}
for all $x\in F := \big(\bigcap_{k=1}^n\dom_{i_k}\big) \cap \big(\bigcap_{\ell=1}^m\dom_{j_l}\big)$, then
$\sum_{k=1}^n a_k\cdot\mu(i_k)\leq \sum_{\ell=1}^m b_\ell\cdot\mu(j_\ell)$.
\item The assignment-routine $\int\_ d\mu \colon S(I, \B \Lambda(X)) \sto \Real$, defined by
\[ (a_k,i_k)_{k=1}^n \mapsto \sum_{k=1}^n a_k\cdot\mu(i_k),\]
is a function.
\end{enumerate}
\end{lem}
\noindent
The next lemma is a slight improvement of Lemma 10.8 in~\cite{BB85}
that will allow us to proceed without using the axiom of countable
choice in the proof of Theorem~\ref{thm: pre-10.10}. It is at this
point that we use the induction principle for disjoint simple
functions as described in Remark~\ref{rem: simpleInduction}.\footnote{We owe this alternative proof to
  a note of the late Erik Palmgren found in the copy of the book~\cite{BB85} by
  Bishop and Bridges that Erik used to own.}

\begin{lem}\label{lem: pre-10.8}
  If $S^+(I, \B \Lambda(X)):=\{v\in S(I, \B \Lambda(X)) \mid \F f_v\geq0\}$ is the set of positive
  simple functions, then there is a function
  ${\phi \colon \mathbb{N}^+ \times S^+(I, \B \Lambda(X))\rightarrow I}$, with
  $(N, v) \mapsto \phi_N(v) \in I$, such that for every $N \in \Nat^+$ and every $v \in S^+(I, \B \Lambda(X))$,
  the following conditions hold:
  \begin{enumerate}[label=(\roman*)]
  \item $\dom_{\phi_N(v)} \subseteq \dom_0(v)$.
  \item $\forall_{x\in \lambda_0^0(\phi_N(v))}\big(\F f_v(x)<N^{-1}\big)$.
  \item $\mu(\phi_N(v))\leq 2N\int v\;d\mu$.
  \end{enumerate}
\end{lem}

\begin{proof}
  Let $N\in\Nat^+$. For the base case let $a\in\Real_{\geq 0}$ and $i\in I$. We
  construct $\phi_N(a,i)\in I$ satisfying
  \begin{itemize}
  \item $\lambda_0^1(\phi_N(a,i))\cup\lambda_0^0(\phi_N(a,i))\subseteq\lambda_0^1(i)\cup\lambda_0^0(i)$,
  \item $\forall_{x \in \lambda_0^0(\phi_N(a,i))}\big(a\chi_{i}(x)<N^{-1}\big)$,
  \item $\mu(\phi_N(a,i))\leq 2Na\mu(i)$.
  \end{itemize}
  Since $(2N)^{-1}<N^{-1}$ we get that $a<N^{-1}$ or $a>(2N)^{-1}$ and using an algorithm that lets us
  decide which case obtains (using Corollary 2.17 in~\cite{BB85}) we set
  \begin{align*}
    \phi_N(a,i)\; := \;
    \begin{cases}
      i\sim i,\;&\text{if}\; a<N^{-1} \\
      i,\;&\text{if}\; a>(2N)^{-1}.
    \end{cases}
  \end{align*}
  Using the fact that $\B\lambda_0(i\sim i)=(\emptyset,\lambda_0^1(i)\cup\lambda_0^0(i))$,
  the verifications of the above properties become routine for both possible values
  of $\phi_N(a,i)$. For the inductive step we assume that we have a disjoint
  $v=(a_k,i_k)_{k=1}^n\in S^+(I, \B \Lambda(X))$ satisfying the
  above conditions and $a_{n+1}\geq 0$ and $i_{n+1}\in I$ disjoint from all the $i_1,...,i_n$.
  Let $w=(a_k,i_k)_{k=1}^{n+1}$ Working similarly, we set
  \begin{align*}
    \phi_N(w)\; := \;
    \begin{cases}
      \phi_N(v) \vee (i_{n+1}\sim i_{n+1}),\;&\text{if}\;a_{n+1}<N^{-1} \\
      \phi_N(v) \vee i_{n+1},\;&\text{if}\;a_{n+1}>(2N)^{-1}.
    \end{cases}
  \end{align*}
  First, assume that $a_{n+1}<N^{-1}$. We have that
  \begin{align*}
    \B\lambda_0(\phi_N(w)) \;:=\; \big( \lambda_0^1(\phi_N(v)) \cap (\lambda_0^1(i_{n+1})\cup\lambda_0^0(i_{n+1}))
    \;,\; \lambda_0^0(\phi_N(v)) \cap (\lambda_0^1(i_{n+1})\cup\lambda_0^0(i_{n+1})) \big).
  \end{align*}
  The first condition then follows immediately from the inductive hypothesis.
  Now let $x\in\phi_N(w)$ and observe that this means that either
  $x\in\lambda_0^0(\phi_N(v)) \cup \lambda_0^1(i_{n+1})$ or
  $x\in\lambda_0^0(\phi_N(v)) \cup \lambda_0^0(i_{n+1})$. In the first case we get
  that $\F f_w(x)=a_{n+1}<N^{-1}$ by our disjointness assumption and in the second case
  we get $\F f_w(x)= \F f_v(x)<N^{-1}$ by the inductive hypothesis. Finally, we get
  \begin{align*}
    \mu(\phi_N(w))\leq \mu(\phi_N(v)) + \mu(i_{n+1} \sim i_{n+1}) \leq \mu(\phi_N(v))\leq
    2N\int v\;d\mu \leq 2N\int w\;d\mu.
  \end{align*}
  Next we assume that $a_{n+1}>(2N)^{-1}$. We have that
  \begin{align*}
    \lambda_0^1(\phi_N(w)) \;&:=\; (\lambda_0^1(\phi_N(v)) \cap \lambda_0^1(i_{n+1}))
    \cup (\lambda_0^1(\phi_N(v)) \cap \lambda_0^0(i_{n+1}))
    \cup (\lambda_0^0(\phi_N(v)) \cap \lambda_0^1(i_{n+1}))  \\
    \lambda_0^0(\phi_N(w)) \;&:=\;\lambda_0^0(\phi_N(v)) \cap \lambda_0^0(i_{n+1})
  \end{align*}
  By the inductive hypothesis the first two conditions follow easily, and for
  the third we get
  \begin{align*}
    \mu(\phi_N(w))\leq \mu(\phi_N(v)) + \mu(i_{n+1}) \leq 2N\int v\;d\mu + 2Na_{n+1}\mu(i_{n+1})=
    2N\int w\;d\mu.
  \end{align*}
 It is straightforward to check that the properties of cases (i)-(iii) are extensional.
\end{proof}

\noindent 
The last lemma needed for the proof of Theorem~\ref{thm: pre-10.10} corresponds to Lemma 10.9 in~\cite{BB85},
and although it reads similar to Lemma~\ref{lem: pre-10.8}, the way it is used
in the proof of Theorem~\ref{thm: pre-10.10} does not invoke countable choice
and it doesn't allow for an induction proof. Hence, we state it without proof,
as a rather direct translation of Lemma 10.9 in~\cite{BB85}.

\begin{lem}\label{lem: pre-10.9}
Let $v:=(a_k,i_k)_{k=1}^n\in S(I, \B \Lambda(X))$ and $c>0$, such that
$\F f_v\leq c$ on $\dom_0(v)$. If
$i\in I$, such that $\F f_v\leq 0$ on $\lambda_0^0(i)\cap\dom_0(v)$, then for every $\varepsilon>0$
there is $j\in I$ satisfying the following conditions:
\begin{enumerate}[label=(\roman*)]
\item $\dom_j \subseteq\dom_0(v)$.
\item $\forall_{x\in\lambda_1^0(j)}\big(\F f_v(x)>\varepsilon\big)$.
\item $\mu(j)\geq c^{-1}\big(\int v\;d\mu - 2\varepsilon\mu(i)\big)$.
\end{enumerate}
\end{lem}
\noindent
Putting everything together we can prove the main result of this section.
The proof follows closely \cite{BB85} but avoids countable choice by using our
Lemma~\ref{lem: pre-10.8} instead of Lemma 10.8 of \cite{BB85} at
the corresponding point in the proof. We refer to~\cite{Ze19} for details.

\begin{thm}\label{thm: pre-10.10} The structure
  $\big(X,S(I, \B \Lambda(X)), \Simple(\B \Lambda(X)), \int\_d\mu\big)$ is a pre-integration space.
\end{thm}

\section{Canonically integrable functions}
\label{sec: L1}

One of the most central constructions in $\BCMT$ is the completion or,
to be more precise, the $L^1$-completion of a Bishop-Cheng integration
space.  To avoid the impredicativities of this definition within
$\BCMT$, we first present in this section the \emph{canonically
integrable functions} explicitly as a family of partial functions. We
then show that this family admits the structure of a pre-integration
space and explain in what sense it can be seen as the completion of
our original pre-integration space.  We follow closely Section 2 of
Chapter 6 in \cite{BB85}, with the exception that we make almost no
mention of \emph{full sets}.  This is because quantification over full
sets is not allowed, even though the property of being a full set can
be defined predicatively.  We will discuss this in more detail below.
As a result, a few of the key lemmas in~\cite{BB85} are missing in our
setting, making some of the proofs, like the one of Theorem~\ref{thm:
  Leb}, more tedious.  We start by giving some basic results on
pre-integration spaces, which we will only state without proof, as
those work completely analogous to the ones in~\cite{BB85},
pp.~217--218.
\textit{For the remainder of this section we fix a pre-integration space $\C L_0 := \big(X, I, \Lambda(X,\Real), \int\big)$}.

\begin{lem}\label{lem: pre-1.2}
\hfill
\begin{enumerate}[label=(\roman*)]
\item Let $i\in I$ and $\alpha \colon \mathbb{N}^+\to I$, such that for
all $n\in\mathbb{N}^+$ we have $\F f_{\alpha_n}\geq 0$ and
$\sum_{n \in \Nat^+}\int\alpha_n \in \Real$, and $\int i+\sum_{n \in \Nat^+}\int\alpha_n>0$. Then
there exists $x\in \lambda_0(i) \cap \bigcap_{n \in \Nat^+}\lambda_0(\alpha_n)$ such that
$$\sum_{n \in \Nat^+} \F f_{\alpha_n}(x) \in \Real \ \ \ \& \ \ \ \F f_i(x)+\sum_{n \in \Nat^+} \F f_{\alpha_n}(x)>0.$$
\item $\forall_{i\in I}\big(\F f_i\geq 0\;\Rightarrow\;\int i\geq 0\big)$.
\item $\forall_{i\in I}\big(|\int i|\leq \int |i|\big)$.
\item If $i,j\in I$ such that $\F f_i(x)\leq \F f_j(x)$, for every $x\in \lambda_0(i) \cap \lambda_0(j)$, then $\int i\leq\int j$.
\end{enumerate}
\end{lem}
\noindent
In classical measure theory, two functions in
$L^1$ are identified, if they agree almost everywhere.
In $\BCMT$, two integrable functions in the $L^1$-completion
of an integration space are identified, if they agree on a full set.
In~\cite{BB85} each function $\B f$ in $L^1$ comes with a representing sequence
$(\B f_n)_n$ of functions from the base integration space.\footnote{This approach to the definition of $L^1$
was developed by Bishop and Cheng in~\cite{BC72} a few years prior to Mikusi\'nski's similar approach to $L^1$ within the classical Daniell integration theory (see~\cite{Mi78, Mi89}).}
Each representing sequence defines the \emph{canonically integrable function} $\sum_n \B f_n$ on a full domain,
and the represented function $\B f$ agrees with $\sum_n \B f_n$
on this domain, i.e.\ they are identified in $L^1$.
Classically speaking, each equivalence class of $L^1$ contains
a canonically integrable function given by
the representing sequence of an element of the equivalence class.
This means that without loss of generality, we can describe $L^1$ predicatively by focusing only on representing sequences and
their associated canonically integrable functions.

\begin{defi}\label{def: canonical}
The set of representations of $\C L_0$ is the totality
\begin{align*}
I_1:=\bigg\{\alpha \in \D F(\mathbb{N}^+, I) \ \mid \  \sum_{n=1}^\infty\int|\alpha_n| \in \Real\bigg\},
\end{align*}
Let $\nu_0 \colon I_1\sto \D V_0$ be given by
\begin{align*}
  \nu_0(\alpha):=\bigg\{x \in \bigcap_{n=1}^\infty\lambda_0\big(\alpha_n\big) \ \mid \ \sum_{n=1}^\infty|\F f_{\alpha_n}(x)| \in\Real \bigg\}.
\end{align*}
Furthermore, let $\mathcal{H} \colon \bigcurlywedge_{\alpha\in I_1}\mathbb{F}(\nu_0(\alpha),X)$
be the dependent assignment routine where $\C H_\alpha:\nu_0(\alpha)\hookrightarrow X$ is the
canonical embedding induced by the embeddings $\lambda_0(\alpha_n)\hookrightarrow X$ and
the dependent assignment routine
$\F g \colon \bigcurlywedge_{\alpha\in I_1}\mathbb{F}(\nu_0(\alpha),\Real)$
given by $\F g_\alpha(x) :=\sum_{n=1}^\infty \F f_{\alpha_n}(x)$, for every $x\in\nu_0(\alpha)$.
We now take $I_1$ to be equipped with the equality
\begin{align*}
  \alpha =_{I_1}\beta :\TOT \big(\nu_0(\alpha), \C H_\alpha, \F g_{\alpha}\big) =_{\F F^{\se}(X)} \big(\nu_0(\beta), \C H_\beta, \F g_{\beta}\big)
\end{align*}
and define the set of \emph{canonically integrable functions}
as the  $I_1$-set of strongly extensional\footnote{By Proposition~\ref{prp: sextpf1}.},
partial functions $\bm\Lambda_1:=(\nu_0,\mathcal{H},\nu_1,\F g)$, where $\nu_1$ is defined through
dependent unique choice as explained in Definition~\ref{def: setofpartialfunctions}.

Let the canonical embedding of $I$ into $I_1$ be the assignment routine $h \colon I \sto I_1$,
defined by the rule $i \mapsto \big(i,0\cdot i,0\cdot i, \ldots \big)$.
\end{defi}

Clearly the assignment routine $h$ is an embedding, since
\begin{align*}
i=_I j & :\TOT \B f_i =_{\F F^{\se}(X)} \B f_j \\
& \Leftrightarrow \bigg(\lambda_0(i) ,\C E_i, \F f_i+\sum_{n=2}^\infty 0\cdot \F f_i\bigg)
=_{\F F^{\se}(X)} \bigg(\lambda_0(j), \C E_j, \F f_j + \sum_{n=2}^\infty 0\cdot \F f_j\bigg) \\
& \Leftrightarrow h(i)=_{I_1} h(j),
\end{align*}
as one can easily verify that $\lambda_0(i) \subseteq \nu_0(h(i))$ and $\lambda_0(j) \subseteq \nu_0(h(j))$.

Following~\cite{BB85}, p.~224, and with a bit of abuse of notation,
we can define basic functions on $I_1$ such as
\begin{align*}
  \_+\_ \colon I_1\times I_1\rightarrow I_1, \ \ \ \ \alpha +\beta :=(\alpha_1, \beta_1, \alpha_2, \beta_2, \ldots),
\end{align*}
satisfying
\begin{align*}
  \F g_{\alpha +\beta} =_{\F F^{\se}(X)} \F g_\alpha + \F g_\beta \ \ \ \ \& \ \ \ \ h(i+j)=_{I_1}h(i) + h(j).
\end{align*}
Similarly, we obtain functions $\_\cdot\_ \colon \mathbb{R}\times I_1 \rightarrow I_1$ and
$|\_ |,~\wedge_1 \colon I_1 \rightarrow I_1$ commuting with their counterparts on $I$ and
the corresponding operations on $\F F^{\se}(X)$.
Note that for construction of these sequences no choice principles are needed.
Finally, the integral $\int \colon I_1 \to \mathbb{R}$ is given by
$$\int \alpha :=\sum_n\int\alpha_n,$$
It is clear that $\int h(i)=\int i$ for all $i\in I$,
which justifies our overloaded notation

The proof of the next lemma follows section 2 of chapter 6 in~\cite{BB85}.

\begin{lem}\label{lem: pre-2.3}
\hfill
\begin{enumerate}[label=(\roman*)]
\item $\forall_{\alpha\in I_1}\big(\big|\int\alpha\big| \leq \int |\alpha|\big)$.
\item If $\alpha\in I_1$, such that $\forall_{x \in \nu_0(\alpha)}\big(\F g_\alpha(x)\geq 0\big)$, then $\int\alpha\geq 0$.
\item If $\alpha,\beta\in I_1$, such that $\forall_{x\in\nu_0(\alpha)\cap\nu_0(\beta)}\big(\F g_\alpha(x)\leq \F g_\beta(x)\big)$, then $\int \alpha\leq\int\beta$.
\item There is a function $\psi \colon I_1 \times \mathbb{N}^+ \rightarrow I_1$, such that for every $\alpha\in I_1$ and $n\in\mathbb{N}^+$, $\psi(\alpha,n) =_{I_1}\alpha$ and
  $$\sum_{k \in \Nat^+}\int |\psi(\alpha,n)_k| \leq 2^{-n} + \int |\alpha|.$$
\end{enumerate}
\end{lem}
\noindent
Lemma~\ref{lem: pre-2.3}(iv) is formulated in a way that allows us to avoid countable choice,
by explicitly constructing function $\psi$.
Unlike in the previous section, we can however still follow
the proof of Lemma 2.14 in~\cite{BB85}.
We are now able to prove the predicative version of
\textit{Lebesgue's series theorem}. The proof generally follows the proof of
Theorem 2.15 in~\cite{BB85}, but we have to be a bit more cautious,
since we don't have a set of a full sets at hand.
For a subset $A$ we can predicatively define what it means to be full,
namely $\exists_{\alpha\in I_1}\nu_0(\alpha)\subseteq A$. However, the totality of full sets is still
defined through separation from $\C P(X)$ and quantification over full sets is thus not possible.

\begin{thm}\label{thm: Leb}
Let $\Gamma \colon \mathbb{N}^+ \to I_1$, such that $\sum_{n \in \Nat^+}\int |\Gamma_n| \in \Real$, and
\[ A:=\bigg\{\; x\in\bigcap_{n=1}^\infty\nu_0\big(\Gamma_n\big) \ \mid \
\sum_{n=1}^\infty |\F g_{\Gamma_n}(x)| \in \Real\bigg\} \]
Then there exists $\alpha\in I_1$ such that $\nu_0(\alpha) \subseteq A$ (i.e.\ A is full) and
\[\forall_{x\in\nu_0(\alpha)}\bigg(\F g_\alpha(x)=\sum_{n=1}^\infty \F g_{\Gamma_n}(x)\bigg).\]
Moreover, if $\alpha\in I_1$ fulfills the above condition, then
$ \lim_{N\rightarrow\infty}\int \big|\alpha-\sum_{n=1}^N\Gamma_n\big| =0$.
\end{thm}

\begin{proof}
We only give a proof sketch and refer the reader to the proof of
Theorem 4.3.12 in~\cite{Ze19} for details.
For each $n\in\mathbb{N}^+$ let $\beta_n:=\psi(\Gamma_n,n)$ with $\psi$
as in Lemma~\ref{lem: pre-2.3}(iv) i.e., $\beta_n\in I_1$, such that for all
$n\in \mathbb{N}^+$ we have $\beta_n=_{I_1}\Gamma_n$ and
$$\sum_{k=1}^\infty\int |\beta_{nk}| < 2^{-n} + \int |\Gamma_n|.$$
It follows that
$\sum_{n=1}^\infty\sum_{k=1}^\infty\int |f_{\beta_{nk}}|
\in \Real$. Let
\begin{align*}
  B:=\bigg\{x\in\bigcap_{n\in\mathbb{N}^+}\bigcap_{k\in\mathbb{N}}\lambda_0(\beta_{nk}) \ \mid \
\sum_{n=1}^\infty\sum_{k=1}^\infty |f_{\beta_{nk}}(x)| \in \Real\bigg\}
\end{align*}
and fix a suitable bijection
$\varphi \colon \mathbb{N}^+ \rightarrow\mathbb{N}^+ \times\mathbb{N}^+$ (e.g.\ as in section 2.3 of~\cite{Ze19}).
Let $\alpha :\mathbb{N}^+ \to I$ be given by
$\alpha_n:=\beta_{\pr_1(\varphi(n))\;\pr_2(\varphi(n))}$,
then\footnote{If $(x_{nk})_{n, k \in \Nat^+}$ is a sequence of sequences of reals and if
  $(y_m := x_{\pr_1(\varphi(m)) \pr_2(\varphi(m))})_{m \in \Nat^+}$, then $\sum_n \sum_k x_{nk}$ converges absolutely
  if and only if $\sum_m y_m$ converges absolutely, and in this case the two sums are equal.
  This fact can be proven constructively and without choice principles for a
  concrete, suitably chosen $\varphi$, see Lemma 2.3.2 in~\cite{Ze19}.}
$\sum_{n=1}^\infty\int |\alpha_n|=\sum_{n=1}^\infty\sum_{k=1}^\infty\int |\beta_{nk}| \in \Real$
and hence $\alpha\in I_1$.
Using the same argument about double series, we can construct an equality of partial functions:
\begin{align*}
\Big(\nu_0(\alpha),~\C H_\alpha,~\F g_\alpha\Big) =_{\F F^\se(X)}\Big(B,~i_B,~\sum_{n=1}^\infty\sum_{k=1}^\infty \F f_{\beta_{nk}}\Big)
\end{align*}
The moduli of equality
$\nu_1(\beta_n,\Gamma_n)$ give inclusions $\nu_0(\beta_n) \eto \nu_0(\Gamma_n)$ for $n\in\Nat^+$
and induce an embedding $e \colon B \eto A$ such that the following diagram commutes
\[
\begin{tikzcd}
&&X \\
\nu_0(\alpha) \arrow[r,equal]\arrow[rdd, bend right, "\F g_\alpha=\sum_n
\F f_{\alpha_n}"'] &B\arrow[rr, "e", hook]\arrow[ru, "i_B", hook]\arrow[dd, "\sum_n\sum_k \F
f_{\beta_{nk}}"] &&A\arrow[lu, "i_A"', left hook->]\arrow[lldd, bend left, "\sum_n \F g_{\,\Gamma_n}"]
\\
\\
&\mathbb{R}.
\end{tikzcd}
\]

To show the second part of the theorem, let $N\in\mathbb{N}^+$ and $\alpha\in I_1$ such that
$\alpha$ satisfies the conditions of the first part of the theorem and set
$\gamma := \big|\alpha-\sum_{n=1}^N\Gamma_n\big| \in I_1$.
If $\delta \colon \mathbb{N}^+ \to I$ is an enumeration of the terms
\[
\begin{matrix}
-\gamma_1 &-\gamma_2 &-\gamma_3 & \cdots \\
|\beta_{(N+1)\;1}| &|\beta_{(N+1)\;2}| &|\beta_{(N+1)\;3}| &\cdots \\
|\beta_{(N+2)\;1}| &|\beta_{(N+2)\;2}| &|\beta_{(N+2)\;3}| &\cdots \\
\vdots &\vdots &\vdots &\ddots
\end{matrix}
\]
into a single sequence using the bijection
$\varphi \colon\mathbb{N}^+ \rightarrow\mathbb{N}^+ \times\mathbb{N}^+$, then
$\sum_{n=1}^\infty\int |\delta_n| =\sum_{n=1}^\infty\int |\gamma_n| \ + \ \sum_{n=N+1}^\infty\sum_{k=1}^\infty\int |\beta_{nk}| \in \Real$, i.e.\ $\delta\in I_1$.
Following the proof in~\cite{BB85} (p.\ 229), for each $x\in \nu_0(\delta)$ we get that
\begin{align*}
\sum_{n=1}^\infty \F f_{\delta_n}(x)=\sum_{n=N+1}^\infty\sum_{k=1}^\infty |\F f_{\beta_{nk}}|(x)-
\sum_{m=1}^\infty \F f_{\gamma_m}(x) \geq 0
\end{align*}
By Lemma~\ref{lem: pre-2.3}(ii) it follows that $\int\delta\geq 0$. Hence
\begin{align*}
0&\leq \int \bigg(\bigg|\alpha-\sum_{n=1}^N \Gamma_n\bigg|\bigg)=\int\gamma =\sum_{n=1}^\infty\int\gamma_n \\
&\leq\sum_{n=N+1}^\infty\sum_{k=1}^\infty\int |\beta_{nk}|\leq \sum_{n=N+1}^\infty\bigg(2^{-n}
 + \int |\Gamma_n|\bigg),
\end{align*}
and the last expression converges to $0$ for $N\rightarrow\infty$.
\end{proof}

\begin{cor}\label{cor: pre-2.15}
If $\alpha \in I_1$, then $\lim_{N\rightarrow\infty}\int \big|\alpha-\sum_{n=1}^N\alpha_n\big| =0$.
\end{cor}
\noindent 
With Lebesgue's series theorem at hand we can now show that the canonically integrable functions
form a pre-integration space,
and as such the complete extension of the pre-integration space $\C L_0$.
All these proofs follow closely \cite{BB85}  so we will omit them altogether.
The final Theorem 2.18 of section 2 of chapter 6 of \cite{BB85} becomes:

\begin{thm}\label{thm: pre-2.18}
$\big(X,I_1,\bm\Lambda_1,\int\big)$ is a pre-integration space.
\end{thm}

\noindent 
In order to treat $L^1$ as the completion of $\C L_0$, we introduce the $1$-norm of $\C L_0$.
In classical measure theory one often identifies integrable functions
that agree almost everywhere and the normed space $L^1$ is
defined modulo this equivalence relation. The positive, constructive
counterpart of this is to identify functions in the complete extension of an
integration space that agree on a full set. Proposition 2.12 in~\cite{BB85}, p.~227,
then tells us that we can define the $1$-norm modulo
this equality. Since in our predicative setting, we don't have recourse to a set of full
set, we need to introduce the $1$-norm a bit
differently. The following fact is straightforward to show.

\begin{prop}\label{thm: pre-int norm}
  Let $p\in I$, such that $\int p=1$.
  \begin{enumerate}[label=(\roman*)]
  \item If $i, j \in I$, the relation $i=_{\int}j :\TOT \int |i-j|=0$
is an equivalence relation on $I$.
  \item The assignment routine
$\int \colon (I,=_{\int}) \sto \mathbb{R}$, given by the rule $i\mapsto\int i$
is a function.
  \item The functions $\cdot$ and $+$ turn
$(I,=_{\int})$ into an $\mathbb{R}$-vector space with neutral element
$0\cdot p$.
  \item The function $||\_||_1 \colon I \rightarrow \mathbb{R}_{\geq 0}$, given by the rule
$$||\,i\,||_1 :=\int |i|,$$
is a norm on $\big((I,=_{\int})\;,\;\cdot,+,0\cdot p\big)$.
  \end{enumerate}
\end{prop}
\noindent
Putting everything together, and in correspondence to Corollaries
2.16, 2.17 in \cite{BB85}, we get the following.

\begin{thm}\label{thm: pre-2.16}
  \hfill
  \begin{enumerate}[label=(\roman*)]
  \item The canonical embedding $h \colon I \eto I_1$ is norm-preserving.
  \item $(I,=_{\int},||\_||_1)$ is a
dense subspace of $(I_1,=_{\int},||\_||_1)$ through $h$.
  \item $I_1$ is complete with respect to $||\_||_1$.
  \end{enumerate}
\end{thm}

\section{Concluding remarks and future work}
\label{sec: concl}

We presented here the first steps towards a predicative reconstruction
$\PBCMT$ of the original impredicative Bishop-Cheng theory of measure
and integration $\BCMT$. Based on the theory of set-indexed families
of sets within $\BST$, we studied the notions of a pre-measure and
pre-integration space, as predicative reformulations of the notions of
a measure and integration space in $\BCMT$.  As first fundamental
examples we presented
\begin{enumerate}[label=(\roman*)]
\item the Dirac measure as a pre-measure,
\item the pre-integration space associated to a locally
compact metric space with a modulus of local compactness, and
\item the pre-integration space of simple functions generated by a pre-measure space.
\end{enumerate}
Finally, we gave
a predicative treatment of $L^1$ as an appropriate completion of the
pre-integration space of the canonically integrable functions.
Using arguments that avoided the use of
full sets and the principle of countable choice, we managed to prove a predicative
version of the constructive Lebesgue's series theorem.

A predicative definition of $L^1$ ensures that all concepts defined
through quantification over $L^1$ in $\BCMT$ become predicative in
$\PBCMT$. For example, quantification over $L^1$ is used in the
Bishop-Cheng definition of a full set\footnote{ The property of being
  a full set can indeed be defined predicatively by quantification over
  the set $I_1$. However, the totality of full sets is still defined by
  separation from the class of all subsets and thus itself a proper
  class.}
(see~\cite{BB85}, p.~224), a constructive counterpart to the
complement of a null set in classical measure theory, and in the
Bishop-Cheng definition of almost everywhere convergence
(see~\cite{BB85}, p.~265). Our predicative treatment of $L^1$ is the
first, clear indication that the computational content of measure
theory can be grasped by $\PBCMT$.

Many question arise naturally from our current work. In~\cite{BB85},
pp.~232--236, the measure space of an integration space is
constructed. A complemented subset $\B A$ of $X$ is called
\textit{integrable}, if its characteristic function $\B \chi_{\B A}$
is in $L^1$, and the \textit{measure} $\mu(\B A)$ is the integral
$\int \B \chi_{\B A}$. A predicative treatment of the pre-measure
space induced by a pre-integration space is expected to be given by
describing the intersection $M = L^1 \cap \F F^{\se}(X, \D 2)$ as an
appropriate set of complemented subsets. The exact relation between
the pre-measure space of the pre-integration space of a given
pre-measure space with the original pre-measure space needs to be
determined.  And similarly for the pre-integration space of the
pre-measure space of a given pre-integration space.  One must also
investigate, if the expected (pre-)measure space of the
pre-integration space $\big(X,I_1,\bm\Lambda_1,\int_1\big)$ is
complete, in the sense of a predicative reformulation of the
definition of a complete measure space (see~\cite{BB85}, pp.~288-289
and~\cite{Pe20}, p.~209).

The Radon-Nikodym theorem is a core result of classical measure theory,
according to which, under appropriate conditions, measures
can be expressed as integrals
$$
\nu(A) = \int_A f d\mu
$$
with respect to other measures. Following the Daniell approach, Shilov
and Gurevich offer a classical treatment of the Radon-Nikodym theorem
in~\cite{Sh66}. Although Bishop tackled it already in~\cite{Bi67}, he
humbly admitted that his treatment
``follows the classical pattern, except that it is much messier'', partly due to the trade-off
requirement of posing stronger hypotheses.  In the light of $\BCMT$,
Bridges offered an improved and extended constructive
version~\cite{Br77}, which led to the revised, joint account with
Bishop given in~\cite{BB85}.  The definition of the notion of absolute
continuity of one integral over another one, which is central to this
constructive proof of the Radon-Nikodym theorem, is impredicative. It
requires quantification over all integrable sets, and therefore over
the proper class of complemented subsets.  As $L^1$ is here
predicatively defined, a predicative treatment of the constructive
Radon-Nikodym theorem within $\PBCMT$ is expected to be possible.

Bishop and Cheng introduced \textit{profiles} in~\cite{BC72} as an
auxiliary concept in order to address convergence in the class of
integrable functions. The profile theorem expresses positively the
classical fact that an increasing function on the reals can have at
most countably many discontinuities. At the same time, it is
responsible for an abundant supply of integrable sets within
$\BCMT$. It also implies the uncountability of reals, and since there
are countable sheaf models of reals~\cite{Sp06b}, there is no hope of
proving the profile theorem constructively without employing some
choice principle.  A proof of a choice-free version of the profile
theorem was given by Spitters~\cite{Sp06b}, using Coquand's point-free
version of the Stone representation theorem.  The question whether we
can recover the basic applications of the theory of profiles through a
choice-free variation of its basic notions and results within $\PBCMT$
is an important open problem.

\bibliographystyle{alphaurl}
\bibliography{refs}

\end{document}